\documentclass[11pt]{amsart}
\usepackage{amsmath,amssymb} 
\usepackage{latexsym}
\newtheorem{thm}{Theorem}[section]
\newtheorem{cor}[thm]{Corollary}
\newtheorem{conj}[thm]{Conjecture}
\newtheorem{lemma}[thm]{Lemma}
\newtheorem{prop}[thm]{Proposition}
\theoremstyle{definition}
\newtheorem{deff}{Definition}[section]
\theoremstyle{remark}
\newtheorem{rem}{Remark}[section]
\newtheorem{example}{Example}[section]
\numberwithin{equation}{section}

\newcommand{\thmref}[1]{Theorem~\ref{#1}}
\newcommand{\secref}[1]{\S\ref{#1}}
\newcommand{\lemref}[1]{Lemma~\ref{#1}}
\newcommand{\propref}[1]{Proposition~\ref{#1}}

\newcommand{\remref}[1]{Remark~\ref{#1}}
\newcommand{\conjref}[1]{Conjecture~\ref{#1}}
\newcommand{\defref}[1]{Definition~\ref{#1}}

\newcommand{\nc}{\newcommand}
\nc{\bib}{\bibitem}
\nc{\on}{\operatorname}
\nc{\res}{\operatornamewithlimits{Res}}
\nc{\ahat}{\hat A}
\nc{\ul}{\underline}
\nc{\arr}{\rightarrow}
\nc{\al}{\alpha}
\nc{\C}{{\mathbb C}}
\nc{\Cn}{{\mathbb C}^n}
\nc{\Z}{{\mathbb Z}}
\nc{\R}{{\mathbb R}}
\nc{\N}{{\mathbb N}}
\nc{\Q}{{\mathbb Q}}
\nc{\A}{{\mathfrak A}}
\nc{\I}{{\frak I}}
\nc{\eps}{\varepsilon}
\nc{\W}{{\frak W}}
\nc{\M}{\mathfrak{M}}
\nc{\mgp}{\M^G[\Sigma_P]}
\nc{\sym}{\mathrm{Sym}}

\nc{\Hom}{\mathrm{Hom}}
\nc{\tensor}{\otimes}
\nc{\ssum}{\cdot`!\sum}
\nc{\im}{\on{Im}}
\nc{\todd}{\mathrm{Todd}}
\nc{\vol}{\mathrm{vol}}
\nc{\Gs}{\Gamma^*}
\nc{\Tc}{\Theta^V}
\nc{\Ch}{\C[[\hbar]]}
\nc{\smv}{\vec{\sigma}}
\nc{\mgps}{\M^G(\Sigma,P,\smv)}
\nc{\trc}{\mathrm{Tr}_{\mathrm{can}}}
\nc{\trq}{\mathrm{Tr}_q}
\nc{\tr}{\mathrm{Tr}}
\nc{\isom}{\cong}
\nc{\Sun}{SU(n)}
\nc{\sut}{SU(2)}
\nc{\gaminv}{\gamma^{-1}}
\nc{\Ah}{A_\hbar}
\nc{\fm}{C^\infty(M)}
\nc{\fmh}{C^\infty(M)[[\hbar]]}
\nc{\cmoy}{\cdot_\mathrm{Mp}}
\nc{\Am}{A_\mathrm{Mp}}
\nc{\Cq}{\C_1(q)}
\nc{\cim}{C^\infty(M)}
\nc{\Se}{\Sigma^\varepsilon}
\nc{\lstart}{L_\mathrm{s}}
\nc{\lend}{L_\mathrm{e}}
\nc{\vstart}{V_\mathrm{s}}
\nc{\vend}{V_\mathrm{e}}
\nc{\mathbold}{\mathbf}
\nc{\SMP}{\Sigma_P}
\nc{\smpe}{\Sigma_P^\varepsilon}
\nc{\Smps}{\Sigma^{P\sigma}}
\nc{\mbp}{P}
\nc{\egam}[1]{E_{\Gamma(#1)}}
\nc{\fgam}[1]{F_{\Gamma(#1)}}
\nc{\vgam}[1]{V_{\Gamma(#1)}}
\nc{\agam}[1]{A_{\Gamma(#1)}}
\nc{\vinc}{V^{\mathrm{inc}}}
\nc{\vout}{V^{\mathrm{out}}}
\nc{\rib}{R}
\nc{\rvsp}{RV}
\nc{\hrib}{HR}
\nc{\func}{Func}
\nc{\cgam}{C_\Gamma}
\nc{\vin}{V_\mathrm{in}}
\nc{\vki}{V_\mathrm{out}}
\nc{\Tr}{\mathrm{Tr}}
\nc{\pin}{\pi^{-1}}
\nc{\fri}{\mathbf{f}}
\nc{\mto}{\mapsto}
\nc{\ue}{\mathfrak{U}}
\nc{\mbc}{\mathbf{c}}
\nc{\mbC}{\mathbf{C}}
\nc{\mC}{\mathbf{CC}}
\nc{\less}{\setminus}
\nc{\lig}{\mathfrak{g}}
\nc{\rg}{\mathrm{Irrep}(G)}
\nc{\rue}{\mathrm{Rep}(\ue)}
\nc{\db}{\partial D}
\nc{\fgsp}{F(G,\smp)}
\nc{\fqgsp}{F^q(\ue,\SMP)}
\nc{\po}[2]{\{#1,#2\}}
\nc{\sign}{\mathrm{sign}}
\nc{\fmgp}{F(\mgp)}
\nc{\fqmgp}{F^q(\mgp)}
\nc{\fb}{\mathbold{f}}
\nc{\gb}{\mathbold{g}}
\nc{\bg}{g}
\nc{\td}{\mathrm{todd}}
\nc{\bc}{\mathbold{c}}
\nc{\Gc}{\check\Gamma}
\nc{\hdual}{\mathfrak{h}^*}
\nc{\n}{\mathfrak{n}}
\nc{\h}{\mathfrak{h}}
\nc{\id}{\mathrm{id}}
\nc{\ut}{\underline t}
\nc{\preg}{\Omega_\mathrm{reg}}
\nc{\ad}{\mathrm{Ad}}
\nc{\bk}{\bar k}
\nc{\qdim}{q\mathrm{dim}}
\nc{\re}{\mathrm{Re}}
\nc{\lt}{\mathfrak{t}}
\nc{\tdual}{\lt^*}
\nc{\dplus}{|\Delta^+|}
\nc{\irr}{\mathrm{Irrep}}
\nc{\orb}{\mathrm{Conj}}
\nc{\holp}{\mathrm{Hol}_p}
\nc{\mon}{\mathrm{Mon}}
\nc{\monp}{\mathrm{Mon}'}
\nc{\rank}{\mathrm{rank}}
\nc{\Smp}{\Sigma\setminus\{p\}}
\nc{\xh}{{\hat{x}}}
\nc{\smp}{\Sigma_p}
\nc{\osum}{\oplus}
\nc{\cyc}{{\mathrm{cycl}}}
\nc{\fel}{\frac12}
\nc{\csi}{\!\!\!\!*}
\nc{\hh}{[[\hbar]]}
\nc{\diff}{\mathrm{Diff}}
\nc{\wt}{\mathrm{wt}}
\nc{\maxx}{{\mathrm{max}}}
\nc{\tras}{\mathrm{Tr}^\mathrm{asym}_\hbar}
\nc{\fini}{$\Box$}
\nc{\ase}{\mathrm{asymp}}
\nc{\frh}{\mathfrak{h}}
\nc{\fra}{\mathfrak{a}}
\nc{\vew}{\vec w}
\nc{\vex}{\vec x}
\nc{\geu}{\kappa}
\nc{\dq}{D_1(q)}
\nc{\ev}{\mathrm{ev}}
\nc{\ndmf}{\mathrm{End}}

\begin{document}

\title[Traces on deformations of moduli spaces]{Trace functionals on non-commutative deformations of moduli
  spaces of flat connections}
\thanks{The research was supported by NSF and NSA grants}
\author{Philippe Roche, Andr\'as Szenes} \address{Massachusetts
  Institute of Technology, Department of Mathematics, Cambridge, MA
  02139, USA} \email{szenes@math.mit.edu} \address{Laboratoire de
  Physique Math\'ematique et Th\'eorique, Universit\'e Montpellier 2,
  34000 Montpellier, France } \email{roche@lpm.univ-montp2.fr}

\maketitle 

\setcounter{section}{-1}
\section{Introduction}

Let $G$ be a compact connected and simply connected Lie group, and
$\Sigma$ be a compact topological Riemann surface with a point $p$
marked on it. One can associate to this data the moduli space of flat
$G$ connections on the punctured Riemann surface $\Sigma$ denoted by
$\M^G=\M^G[\Sigma_p]$. This space decomposes into a union $\M^G =
\cup_\sigma \M^G(\sigma)$ of moduli spaces of flat connections with
the holonomy around the point $p$ fixed in some conjugacy class
$\sigma$ of $G$. For generic $\sigma$ the space $\M^G(\sigma)$ is a
smooth real algebraic orbifold (manifold if $G=\Sun$), and comes
equipped with a canonically defined symplectic form, which induces an
algebraic Poisson structure on it.

In this paper, we study non-commutative deformations of the spaces of
functions on $\M^G$ and $\M^G(\sigma)$.  We bring together the
algebraic approach initiated and developed in the works of Fock, Rosly
\cite{FR}, Alekseev et. al. \cite{AGS}, and Buffenoir, Roche
\cite{BR}, and the theory of formal deformations of symplectic
manifolds, specifically the index theorem of Fedosov and Nest-Tsygan
\cite{DL,Fe,NT}. 

Our first result is a simple construction of a canonically defined
non-commutative algebra $A_q$ depending on a parameter $q$, from which
$F(\M^G)$, the algebraic functions on $\M^G$, may be recovered by
setting $q=1$.  Substituting $q=e^{2\pi i \hbar}$ one obtains an
algebra $A_\hbar$ over the formal power series $\Ch$ which serves as a
formal deformation of $F(\M^G)$.

The central object of the index theorem of Fedosov and Nest-Tsygan is
a cyclic functional $\tr:A_\hbar\rightarrow \C[[\hbar]]$, called the
{\em canonical trace} , which plays the role of the index of an
elliptic operator in this formal theory.

The focus of our work is a conjectural lifting of the canonical trace,
which takes values in formal power series, to a cyclic functional on
$A_q$ taking values in functions of $q$ holomorphic in the unit disc.
The two traces will be related by an asymptotic expansion at $q=1$.

\textsc{Contents of the paper}. In \secref{sec:3series} we recall the
volume formula of Witten and Verlinde's formula and compute the
asymptotics of the $q$-volume series that we introduce. In
\secref{sec:modspaces} we recall the necessary background about the
topology of the spaces $\M^G$ and $\M^G(\sigma)$. We give a short
introduction to the theory of formal deformations in
\secref{sec:formdef}, and explain how a variation of the theory can be
applied to the case of algebraic manifolds. We start the study of the
algebraic functions on the moduli spaces in \secref{sec:graph} by
giving a construction of the Poisson structure and the Poisson trace
using a graphical technique, mostly based on the ideas and
constructions of \cite{FR,AMR,PR}. We quantize this construction in
\secref{sec:ribbon} by introducing a slightly more geometric variant
of the Reshetikhin-Turaev invariants \cite{RT}. The algebra we obtain
is analogous to that in \cite{AGS,BR}, but our construction is
transparent, geometric, and computationally much more efficient.  Our
main result is contained in \secref{sec:proptraces}, where we prove
the asymptotic correspondence of the traces, to which we alluded
above, in the case of $G=\sut$. We review our results and formulate
the main conjecture which served as a motivation for this article in
\secref{sec:conclusion}. In order to avoid crowding the main body of
the paper, and in a hope to make the paper more accessible for the
reader unfamiliar with the theory of quantum groups, we included an
introduction to a relevant part of the subject in an Appendix
\secref{sec:appendix}.

The goal of our work is carrying out the program outlined in
\secref{sec:alg-loc} for the algebras we construct in
\secref{sec:ribbon}. This involves completing the analysis of
\secref{sec:proptraces} for the higher rank groups and proving the
assumptions described in \secref{sec:conclusion}. These will be the
subjects of two follow-up papers. We should note that the idea of
approaching the Verlinde formulas via formal index theory was raised
by Nest and Tsygan.

With these problems out of the way, one could try to approach
\conjref{thm:conj} about the characteristic class of $A_q$. 

\textsc{Acknowledgments} We are greatly indebted to Pavel Etingof for
his help at all stages of this project, and to Boris Tsygan for his
insights and sugessions. Useful discussions with Jean-Michel Bismut,
Victor Guillemin and Mich\`ele Vergne are also gratefully acknowledged.

\tableofcontents

\section{Three series}
\label{sec:3series}

In this section we introduce three families of series associated to
compact Lie groups. They are all related to the topology of the moduli
spaces of flat connections on Riemann surfaces. The first two, the
Witten series and the Verlinde sums are well-known. We will sketch
their geometric significance in the next section. The last one is the
main object of our study, and its connection to the moduli spaces will
be explained in the subsequent parts of the paper.

Before we proceed, however, we need to fix some notation.

\subsection{Lie theory, notation and preliminaries}
\label{sec:lie}

The notions of Lie theory will play a major role in this article. Here
we set the relevant notation.
\begin{itemize}
\item We will assume that $G$ is a compact, simple, connected and
  simply connected Lie group, with complexified Lie algebra
  $\lig$. For most of the paper we also assume that $G$ is simply laced.
\item Fix a maximal torus $T\subset G$ with Lie algebra $\lt$, whose
  complexification is denoted by $\frh$. The pairing between $\lt$ and
  its dual $\lt^*$ will be denoted by $\langle,\rangle$. We will use
  the shorthand $\xh = \exp(x)$ for the exponential mapping, where
  $x\in \lt$ and $\xh\in T$.
\item Let $\Lambda=\exp^{-1}(e)\subset \lt$ be the unit lattice and
  denote by $\Omega$ its integral dual, the weight lattice in $\tdual$.
  The set of roots will be denoted by $\Delta\subset \Omega$, and the
  coroot corresponding to a root $\alpha\in\Delta$ by $\check\alpha$ .
  We write $e_\lambda$ for the group homomorphism $T\rightarrow U(1)$
  corresponding to the weight $\lambda$. Thus we have
  $e_\lambda(\xh)=e^{2\pi i\langle\lambda,x\rangle}$.
\item The Lie algebra $\lig$ has a symmetric bilinear form, which is
  invariant under the adjoint action of $G$. Such a form is unique up
  to multiplication by a constant and induces an inner product on
  $\lig^*$. The normalization of this inner product is usually fixed
  in such a way that the long roots (and coroots) have square length
  2. The form thus normalized is called the {\em basic} inner product
  and will be denoted by $(.,.)$. It induces a linear isomorphism
  between $\lt$ and $\tdual$ denoted by $x\mapsto x^*$.
\item The Weyl group $W_G$ acts on $T$ and $\lt$ effectively.  Assume
  that a Weyl chamber in $\lt$ has been chosen.  This induces a split
  of the roots into positive and negative: $\Delta =
  \Delta^+\cup\Delta^-$ and the choice of the dominant weights $\Omega^+$.
  As usual, $\rho=\frac12\sum_{\alpha\in\Delta^+}\alpha$, and
  $\theta_G$ is the highest root of $G$. The dual Coxeter number $h_G$
  is the integer defined by $h_G=(\theta_G,\rho)+1$.
\item Denote the conjugacy class of an element $t\in T$ by $\sigma_t$. Let
 $T_\mathrm{reg}$ be the set of regular elements in $T$.  We have
  chosen a dominant chamber in $\tdual$, which induces a choice of a
  chamber in $\lt$ and thus a choice of an alcove (a connected
  component) $\mathfrak{a}$ in $T_\mathrm{reg}$. Then by Lie theory we
  have
\[ \orb(G) \isom T/W_G\quad\text{and}\quad \orb_\mathrm{reg}(G)
\isom \mathfrak{a}. \] Denote by $\log$ the local inverse of the
exponential mapping, which maps the only vertex of $\mathfrak{a}$ at
the identity to $0\in\lt$. Then the mapping $\log^*$ identifies
$\mathfrak{a}$ with the following subset of the dominant chamber of
$\tdual$:
\begin{equation}
  \label{eq:alcove}
\mathfrak{a}^* = \{\gamma\in\lt^{*+}|\,(\gamma,\theta_G)<1,\,
(\gamma,\alpha)>0, \,\alpha\in\Delta^+\},  
\end{equation}
where $\theta_G$ is the highest root. 
\item The representation ring $R(G)$ has an integral basis
  $\irr(G)=\{\chi_\lambda\}_{\lambda\in \Omega^+}$, where $\chi_\lambda$ is
  the character of $V_\lambda$, the irreducible representation with
  highest weight $\lambda$.
\item Define a partial ordering on the weights, by saying that
  $\lambda\geq\mu$ if their difference is a sum of positive roots,
  i.e. $\lambda-\mu\in\Z^{\geq0}\Delta^+$. This partial order may be
  extended to $\frh^*$, the set of all complex weights.
\end{itemize}

Other conventions:
\begin{itemize}
\item Given a non-degenerate pairing $Q:V\tensor W\rightarrow\C$
  denote by $\delta_Q$ the \emph{diagonal element} $\sum v^i\tensor
  w_i\in V\tensor W$, where $\{v^i,w_j\}$ are a pair of dual bases of
  $V$ and $W$ correspondingly, i.e. $Q(v^i,w_j) = \delta^i_j$. In
  particular, $\delta(V)\in V^*\tensor V$ is the diagonal elements
  with respect to the canonical pairing between $V^*$ and $V$.
\item Underlining a symbol will mean multiplication by $2\pi i$.
\end{itemize}

\emph{For most of this section we assume that $G$ is simply laced.}
The formulas for the non-simply laced cases require minor
modifications of the ones given.

\subsection{The rational case}
\label{sec:wittseries}

For a positive  integer $g$, consider the function on $\lt$ given by
the series 
\begin{equation}
  \label{eq:h-series}
\widetilde W^G_g(x) = \sum_{\lambda\in
  \Omega^+}\frac{\chi_\lambda(\xh)}{(\dim V_\lambda)^{2g-1}}.
\end{equation}
Clearly, $\widetilde W^G_g(x)$ is a function of $\sigma_\xh$ only.
Up to a normalization, to be discussed below, this is the volume series of
Witten \cite{Wi} (cf. \eqref{eq:witten}).

Recall the Weyl character and Weyl dimension formulas:
\begin{equation}
  \label{eq:weyl}
  \chi_\lambda = \frac1\delta\sum_{w\in W_G} \sign(w)
    e_{w(\lambda+\rho)} 
\quad\text{and}\quad \dim V_\lambda = \prod_{\alpha\in
  \Delta^+}\frac{(\lambda+\rho,\alpha)}{(\rho,\alpha)},
\end{equation}
where $\sign:W_G\rightarrow\pm1$ is the antisymmetric character of
$W_G$ and $\delta=e_\rho\prod_{\alpha\in\Delta^+}(1-e_{-\alpha})$ is
the fundamental antisymmetric character.

Multiplying $\widetilde W^G_g(x)$ by $\delta(\xh)$ and an appropriate
$x$-independent constant one obtains the expression
\begin{equation}
  \label{eq:w-def}
W^G_g(x) = c(g,G)\delta(\xh) \widetilde W^G_g(x)=
(-1)^{(g-1)\dplus}|Z_G|^g  B^G_{2g-1}(x),
\end{equation}
where
\begin{equation}
\label{eq:volume-series}
 B^G_m(x)  =
\sum_{\lambda\in \preg}
\frac{e^{\langle\ul\lambda,x\rangle}}
{\prod_{\alpha\in\Delta^+}(\ul\lambda,\alpha)^m}
\end{equation}
Here $\preg$ is the set of regular, not necessarily dominant weights,
$Z_G$ is the center of $G$ and underlining means multiplication by
$2\pi i$.  The constant $c(g,G)$ is defined by \eqref{eq:w-def}.  The
function \eqref{eq:volume-series} is a multi-dimensional Fourier
series, a higher dimensional analog of the classical Bernoulli
polynomials. Such series were studied in a more general context in
\cite{Sz2} where, in particular, it was proved that $B^G_m(x)$ is
a piecewise polynomial function with rational coefficients.  When
$G=\sut$, up to a normalization factor, one recovers the classical
Bernoulli polynomials
\begin{equation}
  \label{eq:vol-sut}
B_m(x) = \sum_{n\neq 0} \frac{e^{x\ul n}}{\ul n^{m}}.  
\end{equation}
Here and in the rest of this paper, we will omit the subscript or
superscript ``$G$'' when $G=\sut$, if this causes no confusion.  To
compute the coefficients of these polynomials one can apply the
following simple lemma \cite{Sz2}:
\begin{lemma} \label{onedim}
  Let $f$ be a rational function of degree $\leq-2$ on $\C$ and let
  $P_f$ be the set of its poles. Then for each $x$, $0\leq x <1$,
\begin{equation}
\sum_{n\in\Z,\,\ul n\notin P_f} \exp(\ul nx) f(\ul n) = \sum_{p\in
  P_f} \res_{u=p} 
\frac{e^{xu}\,du}{1-e^u} f(u).\label{eq:resrational}
\end{equation}
\end{lemma}
In the case $f(u)=u^{-m}$, $g>1$, one recovers the well-known
formula
\[ B_m(x) = \res_{u=0} \frac{1}{u^{m}}
\frac{e^{\{x\}u}\, du}{(1-e^u)}, \] where $\{\}$ stands for fractional
part of a real number.
\begin{rem}
  The Lemma can be extended to the case $\deg(f)>-2$, by assuming
  $0<x<1$ (cf. \cite{Sz2}).
\end{rem}

\subsection{The trigonometric case}
\label{sec:verlsums}

Denote by $\Omega_\Delta$ the root lattice, which is the integral dual
of the center lattice $\exp^{-1}(Z_G)$.  For positive integers $k,g$
and a dominant weight $\lambda\in \Omega_\Delta^+$, consider the
finite sum
\begin{equation} 
\label{eq:ver-sum}
V^G_{g}(\lambda;k) = \frac1{|W_G|}
\sum_{t\in Z_G[k+h_G]\cap T_\mathrm{reg}}
\chi_\lambda(t)\left(\frac{|Z_G|(k+h_G)^{\mathrm{rank}(G)}}
  {\prod_{\alpha\in\Delta}(1-e_\alpha(t))}\right)^{g-1},
\end{equation} 
where $Z_G[k]=\{t\in T|\, t^k\in Z_G\}$, and $h_G=(\theta_G,\rho)+1$.
This sum, first written down by E. Verlinde \cite{Ve}, turns out to be
an integer valued function, whose dependence on $\lambda$ and $k$ is
again piecewise (quasi-)polynomial.  Note that the
denominator is the Weyl density function which can be written as
$\delta\bar\delta=(-1)^{\vert \Delta_{+}\vert }\delta^2$.  Using this
and taking advantage of the Weyl character formula again, we arrive at
\begin{equation}
  \label{eq:v-e}
V^G_{g}(\lambda;k) =
\left((-1)^{\dplus}|Z_G|(k+h_G)^{\mathrm{rank}(G)}\right)^{g-1}
\sum_{t\in Z_G[k+h_G]\cap T_\mathrm{reg} }
\frac{e_{\lambda+\rho}(t)}{\delta(t)^{2g-1}}.
  \end{equation}

For $G=\sut$ we obtain  
\[V_g(l;k) = i(2(k+2))^{g-1} \sum_{j=1}^{2k+3}
\frac{e^{\ul j(l+1)/(2(k+2))}}
{\left(2\sin(j\pi/(k+2))\right)^{2g-1}},\quad j\neq k+2,
\] 
where $l$ is an even number.  Again, we can use residue techniques to
evaluate this sum:
\begin{lemma}
  Let f(z) be a rational function on $\C$ with a set of poles $P_f,$
  $m$ a positive integer, such that $f(z)\,dz/(z(1-z^{m}))$ is regular
  at 0 and at $\infty$. Then
\[ 
\sum_{z^m=1,\, z\notin P_f} f(z) = \sum_{p\in P_f}
\res_{z=p}\;\frac{m\, dz}{z(1-z^m)}f(z).
\]
\end{lemma}
Applying the lemma to our sum with $f(z) = z^{l+1}(z-1/z)^{1-2g}$, and
$m=2(k+2)$, we obtain
\begin{equation}
  \label{eq:ver-res}
V_g(l,k) = -2(-2(k+2))^g\res_{z=1}\;
\frac{z^{l'}\,dz}{(1-z^{2(k+2)})(z-1/z)^{2g-1}z},
\end{equation}
where $l'=l+1\mod 2k+4$ and $0\leq l'<2k+4$. In principle, we also
need the residue at $z=-1$, but by symmetry it is equal to the one at
$z=1$ ($l'$ is odd) hence the extra factor of $2$.

Finally, we can make the substitution $z=e^{u/(2k+4)}$ to arrive
at
\begin{equation}
  \label{eq:v-p}
V_g(l,k) = 2(-2(k+2))^{g-1}P_{2g-1}(\{(l+1)/2(k+2)\},k+2),
\end{equation}
where
\begin{equation}
  \label{eq:def-p}
P_m(x,k) = \res_{u=0}
\frac{e^{xu}\, du}{(1-u)(e^{u/2k}-e^{-u/2k})^m}.
\end{equation}

\subsection{The $q$-rational case}
\label{sec:qrat}

Now we turn to the main object of our study.  Define the $q$-integers
by $[n]_q=\frac{q^n-q^{-n}}{q-q^{-1}}\in \Z[q,q^{-1}]$, for every
$n\in\Z^{>0}$.  The $q$-dimension of $V_\lambda$ is defined in analogy
with \eqref{eq:weyl} by $$\qdim V_{\lambda}=
\prod_{\alpha\in\Delta^+}\frac{[(\lambda+\rho,\alpha)]_q}
{[(\rho,\alpha)]_q}\in \Z[q,q^{-1}].$$

By replacing the classical dimension by the $q$-dimension, we can
write down the $q$-version of the $\widetilde W^G_g$ series:
$$
\widetilde T_g^G(x;q) = \sum_{\lambda\in \Omega^+}
\frac{\chi_\lambda(\hat{x})}{(\qdim V_\lambda)^{2g-1}}.
$$
The function $\qdim V_\lambda$ as a function of $\lambda$ has the
same symmetry properties with respect to the Weyl group as the usual
dimension, so we can use the same trick as above to arrive at the
analog of the $W^G_g$ series: $T^G_g(x;q) = c(g,G,q)\delta(\xh)\widetilde
T^G_g(x;q)$,
\[ T^G_g(x;q) =  (-1)^{(g-1)\dplus}|Z_G|^g
\sum_{\lambda\in \preg} \frac{e^{\langle\ul\lambda,x\rangle}}
{\left(\prod_{\alpha\in\Delta^+}
    q^{(\alpha,\lambda)}-q^{-(\alpha,\lambda)}\right)^{2g-1}}.
\]
Note that while $c(g,G,q)$ is an analog of $c(g,G)$, we do not have
$c(g,G,1)=c(g,G)$.

When $G=\sut$, the formula reads
\begin{equation}
  \label{eq:tqsut}
  T_g(x;q) =  (-1)^{g-1}2^g\sum_{n\neq 0}
\frac{e^{\ul n x}}{(q^{n}-q^{-n})^{2g-1}}.
\end{equation}
Observe that for $g \geq 1$ this series converges to a holomorphic
function on the the unit disc in the complex $q$-plane.  It is
difficult to evaluate such a sum (although, cf.
\remref{thm:thetafunction}). Instead, we will study the behavior of
$T_g(x;q)$ as $q\rightarrow 1$. More precisely, we will compute the
asymptotics of $T_g(x;e^{\pi i\hbar})$ as $\hbar\rightarrow i0^+$,
i.e. as $\hbar$ approaches 0 along the ray of purely imaginary numbers
in the upper half plane. To this end, consider the form
\begin{equation}
\label{eq:def-w}
w_g(u,x;\hbar) = \frac{e^{\{x\} u}\,du}{1-e^{u}}
\frac{1}{(e^{\hbar u/2}-e^{-\hbar u/2})^{2g-1}}
\end{equation}
in the complex $u$-plane for a fixed $\hbar\in i\R^+$.  Divide the set
of its poles into 3 parts as follows.
\begin{enumerate}
\item $P_1=\{\ul n|\;n\in\Z,\,n\neq0\}$;
\item $P_2=\{ \ul n/{\hbar}|\; n\in\Z,\,n\neq0\}$;
\item $u=0$.
\end{enumerate}
Clearly, we have
\[  
T_g(x;e^{\pi i\hbar}) = (-1)^{g-1}2^g\sum_{p\in P_1}
\res_{u=p}w_g(u,x;\hbar).
\] 
Now consider $\sum_{p\in P_2} \res_{u=p}w_g(u,x,\hbar)$ and assume
$x\notin \Z$. Since $e^{\{x\} u}/(1-e^{u})$ and its first $2g$
derivatives vanish exponentially as $u\rightarrow\pm\infty$ on the
real line and $(e^{\hbar u/2}-e^{-\hbar u/2})^{2g-1}$ is periodic with
period $2\pi/{\hbar}$, we have
$$\left|\res_{u=2\pi n/\hbar} w_g(u,x,\hbar)\right|<ce^{-\tau
  |n|/|\hbar|}$$
for some positive constants $c$ and $\tau$,
independent from $n$. Then by summing the geometric series we obtain
that
$$\left|\sum_{p\in P_2}
  \res_{u=p}w_g(u,x,\hbar)\right|<2ce^{-\tau/|\hbar|}$$
as
$\hbar\rightarrow i0^+$.  Finally, since, outside a small neighborhood
of its poles, the function $e^{\{x\}u}/(1-e^{u})$ vanishes
exponentially as $|\re(u)|\rightarrow\infty$ and so does $(e^{i\hbar
  u}-e^{-i\hbar u})^{1-2g}$ as $|\im(u)|\rightarrow \infty$, we see
that the line integral of $w_g(u,x;\hbar)$ over a sequence of
appropriately chosen contours, e.g. the boundary of the rectangles
\[
\mathrm{Rect}_L = \{|\hbar\re(u)|,|\im(u)|\leq (2L+1)\pi \},\quad L\in
\mathbb{N}
\]
goes to 0 as $L\rightarrow\infty$.  We can summarize what we have
found as follows:
\begin{prop}
  There are positive constants $c$ and $\tau$ (possibly depending on
  $x$) such that for $x\notin \Z$
  \begin{equation}
    \label{eq:qwit-res}
|T_g(x;e^{\pi i\hbar}) - (-1)^{g-1}2^g\res_{u=0}
w_g(u,x;\hbar)|<ce^{-\tau/|\hbar|} 
    \end{equation}
    for sufficiently small $\hbar\in i\R^+$.
\end{prop}

\begin{rem}
\label{thm:thetafunction}
Note that the existence of such an expansion around $q=1$ for a
holomorphic function of $q$ on the unit disc, is a rather rare
occurrence. It strongly suggests that the function is related to
modular forms. In fact, there are examples of such relations
\cite{WW}:
  \[ T_1(x,q) =\frac d{du}\log \theta_4(u/2,\tau), 
  q=e^{i\pi \tau}\] but we will not explore this connection in this
  paper. From the point of view of modular forms, the aymptotic
  behavior of $T$ appears as a ``defect'' of sorts, a measurement of
  its failure to be modular.
\end{rem}

We finish this section with an observation which will be central to
our main result. Clearly, there is a formal analogy between
$V_g^G(\lambda;k)$ and $T_g^G(x;\hbar)$, both of them being a
trigonometric deformation of $W_g^G(x)$. Our residue calculations
quantify this analogy in the case of $G=\sut$, as follows.

First, by shifting the variables we can rewrite \eqref{eq:v-p} as
\begin{equation}
  \label{eq:sh-vp}
V_g(l-1;k-2)=2(-2k)^{g-1}P_{2g-1}({l/2k},k),  
\end{equation}
for $l/2,k\in\Z^{\geq0}$. On the other hand,
\eqref{eq:qwit-res} implies that asymptotically 
\begin{equation}
  \label{eq:t-p}
T_g(x;e^{\pi
  i\hbar})\sim 2(-2)^{g-1}P_{2g-1}(\{x\},1/\hbar)
\quad\text{as}\quad\hbar\rightarrow i0^+.  
\end{equation}

A similar equality holds in the higher rank case. This will be covered
in a later publication.

\subsection{Several punctures}
\label{sec:sev-pun}

Let $P$ be a finite set. We can extend the results of this section as
follows.  We can write down a function of $\vex:P\rightarrow\lt$
\begin{equation}
 \tag{\ref{eq:h-series}P}
 \widetilde W^G_g(\vex) = \sum_{\lambda\in
  \Omega^+}\frac{\prod_{p\in P}\chi_\lambda(\hat{\vec{x}}(p))}{(\dim
  V_\lambda)^{2g-2+|P|}}. 
\end{equation}
and the series
\begin{multline}
\tag{\ref{eq:volume-series}P}
 W^G_g(\vex)= c_{|P|}(g,G)\left(\prod_{p\in P}
   \delta(\vex(p))\right)\widetilde W^G_g(\vex) 
= \\
\frac{(-1)^{(g-1)\dplus}|Z_G|^g}{|W_G|}
\sum_{\vew:P\rightarrow W_G}\sign(\vew)B^G_{2g-2+|P|}(\vew\cdot\vex), 
\end{multline}
where  $\sign(\vew) = \prod_{p\in
  P}\sign(\vew(p))$ and $\vew\cdot\vex = \sum_{p\in
  P}\vew(p)(\vex(p))$.  For $G=\sut$ we obtain the formula
\begin{equation}
\tag{\ref{eq:vol-sut}P}
 W_g(\vex) = (-2)^{g-1}\sum_{\vec\epsilon}\sign(\vec\epsilon)
B_{2g-2+|P|}\left(\vec\epsilon\cdot\vex\right)
\end{equation}
where the first summation is over all possible $|P|$-tuples of signs
$\vec\epsilon$ and $\sign(\vec\epsilon)=\prod_p\vec\epsilon(p)$.

In the trigonometric case, one replaces $\lambda\in \Omega_\Delta$ by
$\vec\lambda: P\rightarrow \Omega^+$ such that $\sum_{p\in P}
\vec\lambda(p)\in \Omega_\Delta$, and writes a formula for
$V_g(\vec\lambda;k)$ by replacing $\chi_\lambda$ in \eqref{eq:ver-sum}
by $\prod_{p\in P}\chi_{\vec\lambda(p)}(t)$. Now the formulas for $\widetilde
T^G_g(\vex;q)$ and $T^G_g(\vex;q)$ can be written down by
analogy. We can summarize the final result as follows:
\begin{prop}
  \label{thm:fin-pun}
  The Verlinde sums and the $q$-volume series for $G=\sut$ are both
  related to the same polynomial function as follows:
  \begin{equation}
\tag{\ref{eq:sh-vp}P}
V_g(\vec l-\vec 1,k-2) = 2(-2k)^{g-1} \sum_{\vec\epsilon}
\sign(\vec\epsilon)P_{2g-2+|P|}\left(\{\vec\epsilon\cdot\vec
  l/2k\},k\right);
\end{equation} 
and if for every choice of signs $\vec\epsilon$ the condition
$\vec\epsilon\cdot\vex\notin\Z$ holds then
\begin{equation}
\tag{\ref{eq:t-p}P}
T_g(\vex,e^{i\pi\hbar}) \sim 2(-2)^{g-1}\sum_{\vec\epsilon}
   \sign(\vec\epsilon) P_{2g-2+|P|}\left(\{\vec\epsilon\cdot\vec
      x\},\hbar^{-1}\right),
\end{equation}
where $\sim$ means asymptotic equality as used above.
\end{prop}
We will call $\hat{\vex}\in T$ \emph{special} if the condition in the proposition 
holds.

\section{Moduli spaces}
\label{sec:modspaces}

This section serves as a quick introduction to the topology of the
moduli spaces of flat connections on Riemann surfaces. For more
detailed analysis \cite{BL} is a good reference. The first part
\secref{sec:topmod} is not essential for following the rest of the
paper and is only given as an orientation for the reader. We want to
emphasize the relation between the the singularities of the Witten
sums and those of the corresponding moduli spaces.

Keeping the notation of the previous section, let again $G$ be a
compact, simple, simply connected Lie group. Let $\Sigma$ be a
topological Riemann surface and $P\subset \Sigma$ be a finite nonempty
set of points. We will use the shorthand $\SMP=\Sigma\setminus P$ and
when $P=\{p\}$ we will write $\smp$ for $\Sigma\setminus P$.  The
moduli space of flat $G$-connections $\M^G=\M^G[\SMP]$ modulo gauge
transformations can be defined as the quotient
\[ \M^G = \Hom(\pi_1(\SMP),G)/\ad G.\]
We discuss the case of one puncture $P=\{p\}$ first, and
explain how to generalize the results to several punctures at the end
of the section.

There is a natural map $\holp:\M^G[\smp]\rightarrow\orb(G)$ which
assigns to a flat connection the conjugacy class of the holonomy
around the puncture $p$.  Denote by $\M^G(\sigma)$ the space
$\holp^{-1}(\sigma)$ for and conjugacy class $\sigma\in\orb(G)$. We will
be interested in the case of regular orbits only. In this case,
$\M^G(\sigma)$ is (a possibly singular) manifold of dimension
$(2g-1)\dim G-\rank(G)$.

\subsection{Topology and singularities}
\label{sec:topmod}

In this paragraph, we describe the family of spaces $\M^G(t)$ as $t$
varies, with emphasis on the question of smoothness. This is
needed since we will discuss the Riemann-Roch calculus on these spaces
later. The discussion is necessarily incomplete and somewhat informal.
For complete details cf. \cite{BL,MW}.  

Choosing the standard presentation of $\pi_1(\smp)$ we can represent
$\M^G(\sigma)$ as
\[ \{[A_1,B_1]\dots[A_g,B_g]\in\sigma|\, A_i,B_i\in G\}/\ad G, \]
where $[A,B]=ABA^{-1}B^{-1}$.  As in the previous section, sometimes
we will replace the regular orbit $\sigma$ by a representative $t\in
T_{\mathrm{reg}}$.  Then we have
\[ \M^G(t) = \{[A_1,B_1]\dots[A_g,B_g]=t\,|\, A_i,B_i\in G\}/\ad T. \]

Since $G$ is compact, the space $\M^G(t)$ is Hausdorff. We will
study the question: at what $\mu\in\M^G$ is the space $\M^G$
singular. For $\mu\in\M^G$, let $\mon(\mu)$ be the monodromy group of
$\mu$, which is the subgroup of $G$ generated by holonomies of $\mu$,
i.e. by the subgroup generated by the elements $\{A_i,B_i|\,
i=1,\dots,g\}$. Assume that $\mon(\mu)$ is connected. Then
$\monp(\mu)$ the commutator subgroup of $\mon(\mu)$ is a semisimple
compact Lie subgroup of $G$ containing $t$, thus it has a maximal
torus $T_\mu\subset T$.  Denote the set of roots of $\monp(\mu)$ by
$\Delta_\mu$ and the Lie algebra of $T_\mu$ by $\lt_\mu$.

Barring some degenerate cases, for a generic $\mu$, we have
$\monp(\mu)=G$. Clearly, the singularities will appear at a solution
$\mu$ whenever the centralizer $Z(\monp(\mu))$ in $G$ is
strictly greater than $Z_G$. Denote $t=\holp(\mu)$. There are two cases:\\
\smallskip

\noindent 1. $\dim Z(\monp(\mu))>0$. In this case we say that $\M^G$
has a {\em serious} singularity at $\mu$. We can assume that $Z(\monp(\mu))$
contains a 1-dimensional toric subgroup $T_1\subset T$. Then
  $$\Delta_\mu\subset\{\alpha\in\Delta|\,e_\alpha(t)=1,\text{ for }
  t\in T_1\}.$$
  Thus
  $$\lt_\mu\subset\langle\check\alpha|\,\alpha\in\Delta_\mu
  \rangle_\mathrm{lin},$$
  where $\langle\rangle_\mathrm{lin}$ means
  linear span and $\check\alpha$ is the coroot corresponding to the
  root $\alpha$.  This happens if and only if $t$ is not in general
  position with respect to the coroots of $G$. In other words,
  $\M^G$ has serious singularities at $\mu$ if $t$ is the exponential of a
  linear combination of $\rank(G)-1$ coroots of $G$. We will call such
  $t$ \emph{special}. This property only depends on the conjugacy
  class $\sigma_t$.\\
  2. $Z(\monp(\mu))$ is finite but strictly larger than $Z_G$. This
  produces orbifold singularities. The existence of such $\mu$ does
  not depend on $t$, but on the group $G$ only.  More precisely, one
  needs a non-trivial element $z\in Z(\monp(\mu))/Z_G$, i.e. an
  element which lies on $\rank(G)$ singular subsets
  $U_\alpha=\{e_\alpha=1\}\subset T$, $\alpha\in\Delta$, but not in
  the intersection of all of the $U_\alpha$s which is $Z_G$.  Note
  that such element $z$ does not exist for $G=\Sun$, thus in this case
  orbifold singularities do not appear in $\M^G$.
  \smallskip
  
  We can conclude that $\M^G$ is singular (in the serious sense) at a
  point $\mu$ exactly when $\holp(\mu)$ is special. Some more work is
  required to show that $\M^G(t)$ is singular exactly when $t$ is
  special \cite[Sec. 5]{BL}. This implies, in particular, that
  $\M^G(t)$ is smooth when $G=\Sun$ and $t$ is not special.

  Recall from \eqref{eq:alcove} that the set of regular orbits in $G$
  is represented by an alcove $\mathfrak{a}\subset T$, which is in
  one-to-one correspondence with the set $\mathfrak{a}^*\subset\tdual$
  of dominant real weights of height at most 1 via the mapping
  $\log^*$.  We can introduce the sets $\mathfrak{a}_\mathrm{nspec}$
  (resp.  $\mathfrak{a}^*_{\mathrm{nspec}}$) as the subsets of $T$
  (resp. $\tdual$) representing nonspecial elements.
  By the above discussion, this set is the complement of the
  intersection of $\mathfrak{a}$ (resp.  $\mathfrak{a}^*$) with a
  hyperplane arrangement, and has a rather complicated chamber
  structure.  The spaces $\M^G(t)$, where $t$ varies in one of the
  chambers of $\mathfrak{a}_\mathrm{nspec}$ are all isomorphic and form
  a trivial family.  The spaces $\M^G(t)$ and $\M^G(t')$ corresponding
  to two neighboring chambers differ by set of high codimension and in
  general are not isomorphic.  For details see \cite{BH,BL,JK,Th}.

\subsection{Line bundles and the symplectic form} 
\label{sec:sympl-form}

A useful approach for studying the topology of $\M^G(t)$ is to represent
it as a quotient using infinite dimensional symplectic reduction
\cite{AB,BL,FR,MW}. In particular, this means that $\M^G(t)$ comes
equipped with a canonical symplectic form $\omega_t$ induced by the
symplectic form on the space of all connections, which depends only on
the normalization of the symmetric bilinear product on $\lig$.  The
line bundles on $\M^G(t)$ correspond to homogeneous bundles for the
central extension of the loop group of $G$. We summarize the necessary
facts in the following proposition. For details consult
\cite{BL,MW}.

\begin{prop} 
\label{thm:top-mod}
Let $t\in\fra_\mathrm{nspec}$ be a non-special element of $T$. Then
  \begin{itemize}
  \item There is an identification $\eta:\tdual\osum\R\isom
    H^2(\M^G(t),\R)$, such that
  \item The cohomology class of $\omega_t$ is $\eta(\log(t)^*,1)$
  \item For $\lambda\in P_\Delta$ and $k\in\Z$, the cohomology class
    $\eta(\lambda,k)$ is integral, being the Chern class of a line
    bundle $L_{\lambda,k}$.
  \item The line bundle $L_{\lambda,k}$ is positive on $\M^G(t)$ with
    respect to $\omega_t$ if and only if $\lambda/k$ is in the same
    chamber as $\log(t)^*$.
  \item The following formula holds:
    $c_1(\M^G(t))=2\eta(\rho,(\theta_G,\rho)+1)$.
  \end{itemize}
\end{prop}
Now we can describe the connection of the formulas of the previous
section with the topology of the moduli spaces. The first formula
proved by Witten \cite{Wi,KL} says that the functions $W^G_g(\log(t))$
give the volume of the moduli spaces with respect to the canonical
symplectic structure.
\begin{equation}
  \label{eq:witten}
  W^G_g(\log(t)) = \int_{\M^G(t)}e^{\omega_t}.
\end{equation}
The second formula computes the Riemann-Roch number of the line bundle
$L_{\lambda,k}$ on the space $\M^G(t)$, where $\lambda/k=\log(t)^*$.
Assume that $(\lambda,k)$ is such that $t$ is non-special. The space
$\M^G(t)$ may be endowed with an appropriate K\"ahler structure and
the line bundle $L_{\lambda,k}$ with a holomorphic structure. Then the
famous Verlinde formula reads (cf. \cite{Be,MW,BL}) 
\begin{equation}
  \label{eq:ver-formula}
 V_g^G(\lambda;k) =  \dim H^0(\M^G(t),L_{\lambda,k}).
\end{equation}
The 4th part of the Proposition corresponds to the statement that the
bundle $L_{\lambda,k}$ is ample on $\M^G(t)$ under the stated
conditions. In particular, vanishing of higher cohomology holds for
sufficiently high powers of $L_{\lambda,k}$. Somewhat surprisingly,
such vanishing holds for the bundles $L_{\lambda,k}$ themselves
\cite{Te}.  Thus by the Grothendieck-Riemann-Roch formula we have
\begin{equation}
  \label{eq:ver-rr}
   V_g^G(\lambda;k) =  \int_{\M^G(t)}e^{\eta(\lambda,k)}\todd(\M^G(t))
\end{equation}
The RHS does not change if we vary $t$ inside its chamber.

Finally, using that $\hat A(M)=e^{c_1(M)/2}\todd(M)$, and comparing
the formulas for $\omega_t$ and $c_1(M)$ given in
\propref{thm:top-mod} with \eqref{eq:sh-vp} we obtain ($G=\sut$):
\begin{equation}
  \label{eq:rr-ahat}
2(-2)^{g-1}P_{2g-1}(x,k) =
\int_{\M^G(\xh)}e^{k\omega_{\xh}}\ahat(\M^G(\xh)),
\end{equation}
for $0<x<1/2$.

A similar formula holds in the for general $G$ as well
\cite{Sz3,BL}.

\subsection{Several punctures}
\label{sec:mod-sev-pun}

The case of several punctures is entirely analogous, so we will be very 
brief. In this case, the moduli space $\M^G[\Sigma_P]$ is a union of
the fibers $\M^G(\vec\sigma)$ of the map
$\mathrm{Hol}_P:\M^G[\Sigma_P]\rightarrow\orb(G)^{\times|P|}$.
Serious singularities arise whenever for some $|P|$-tuple of Weyl
group elements $\vew:P\rightarrow
W_G$ the product $\prod_{p\in P} \vew(\vec t)$ is
a special element of $T$ (as defined above). We will call such
$|P|$-tuple of elements of $T$ $\vec t:P\rightarrow 
\mathfrak{a}$ representing $\vec\sigma$ \emph{special}.
Again, there is a canonical symplectic form $\omega_{\vec t}$ on
$\M^G(\vec t)$ and an isomorphism $\eta:\lt^{\osum|P|}\osum\R\isom
    H^2(\M^G(\vec t),\R)$. The multiple puncture version of the
    formula for $G=\sut$ from above reads:
\begin{equation}
  \tag{\ref{eq:rr-ahat}P}
2(-2)^{g-1} \sum_{\vec\epsilon:P\rightarrow\pm}
\sign(\vec\epsilon)P_{2g-2+|P|}(\{\vec\epsilon\cdot\vex
\},k) = \int_{\M^G(\vec \xh)}e^{k\omega_{\vec \xh}}\ahat(\M^G(\vec \xh)),
\end{equation}
where again $0<\vex(p)<1/2$, for all $p\in P$.

\section{Formal deformations of manifolds and index theorems}
\label{sec:formdef}

\subsection{Formal deformations of symplectic manifolds}

In this section we review the formal deformation theory of symplectic
manifolds. A reference for this is \cite{Fe,We}.

For a complex vector space $V$, denote by $V[[\hbar]]$ the space of
formal power series in $\hbar$ with coefficients in $V$.  A {\em
  formal deformation} (the terms {\em star product} or {\em
  deformation quantization} are also used) of a manifold $M$ is a
product $\cdot:\fm \tensor \fm \arr \fmh$, which when extended
linearly to $\fmh$ is
\begin{itemize}
\item associative: $f\cdot (g\cdot h) = (f\cdot g)\cdot h$;
\item local: $f\cdot g = fg + \sum_{n=1}^\infty B_n(f,g)\ul\hbar^n$,
  where $B_n$ is a bidifferential operator.
\end{itemize}
It is easy to see that $B_1(f,g)-B_1(g,f) = \{f,g\}$ is a Poisson
bracket on $\cim$. When $M=(M^{2n},\omega)$ is a symplectic manifold, and
the deformation is such that the induced Poisson bracket is the
symplectic one, then we speak of a formal deformation of a symplectic
manifold.
\begin{rem}
  Note the factor of $2\pi i$ in front of $\hbar$, marked by
  underlining, inside the expansion in the definition of locality.
  This at variance with most conventions in the literature. We chose
  it because it is consistent with $\hbar$ being real and with the
  symplectic form $\omega$ being integral in our applications.
\end{rem}
Associativity induces an infinite set of complicated non-linear
equations on the $B_n$s, which were recently explicitly solved in
local coordinates by Kontsevich \cite{Ko}. The existence of the
solution in the symplectic case was proved \cite{DL,G} (cf. \cite{We}
for detailed references).

The group of formal base changes (gauge transformations)
\begin{equation}
\mathcal{G}(M) = \left\{\gamma:f\mto f+\sum_{n=1}^\infty
  \gamma_n(f)\left|\right.\; \gamma_n \text{ is a differential
    operator}\right\}  
\label{eq:gauge}
\end{equation}
acts on the space of formal products via $f\cdot_\gamma g =
\gamma(\gaminv f\cdot\gaminv g)$. A natural question is the
classification of formal deformations up to this action.

From here on, we only study the symplectic case $(M^{2n},\omega)$.
Then the orbits of this action are labeled by a {\em characteristic
  class} $\theta(M,\cdot)\in\omega/\hbar + H^2(M)\hh$ (cf.
\cite{DL,NT,KS,WeHo, NT, De}).

Another interesting related object is the \emph{canonical trace} on
the algebra $(\fmh,\cdot)$.  A functional $T$ on a non-commutative
ring is called \emph{cyclic} or a \emph{trace} if $T(a\cdot
b)=T(b\cdot a)$.  Such functionals on $(\fmh,\cdot)$ form a
1-dimensional free module over $\Ch$ with a distinguished element. To
understand this we recall the local theory of formal deformations of
symplectic manifolds.

Let $M=\R^{2n}$ with a translation invariant symplectic form $\omega$.
Then the famous Moyal product is given by
\[ f\cmoy g = m(e^{\ul\hbar\pi_\omega} f\tensor g),
\]
where $\pi_\omega=\omega^{-1}$ is the translation invariant Poisson
bivector field induced by $\omega$ and $m$ is the ordinary commutative
product on $\fmh^{\tensor2}$. We collect the highlights of the local
theory in the theorem below.
\begin{thm}[\cite{Fe,NT}]
  \begin{enumerate}
  \item The Moyal product is a formal deformation of the symplectic
    vector space $(\R^{2n},\omega)$.
  \item Any formal product on $(\R^{2n},\omega)$ is gauge-equivalent
    to $\cmoy$.
  \item For every local derivation $D=\sum_{n=1}^\infty D_n \hbar^n$
    of the algebra $(C^\infty(\R^{2n})\hh,\cmoy)$, $D_n\in
    \diff(\R^{2n})$, there exists an element $f\in
    C^\infty(\R^{2n})\hh$, such that $Dg = f\cmoy g-g\cmoy f$, for
    every $g\in C^\infty(\R^{2n})\hh$. Similarly, any local
    automorphism of the form $1+\sum_{n=1}^\infty A_n \hbar^n$ can be
    obtained by exponentiating such a derivation.
  \item Integration against the density $\omega^n/\hbar^n$ defines a
    trace functional on the algebra $(C^\infty(\R^{2n})\hh,\cmoy)$,
    which is unique up to a constant in $\Ch$.
  \end{enumerate}
\end{thm}
The particular choice of normalization $f \rightarrow \int
f\omega^n/\hbar^n$ is called the \emph{canonical trace} on the Moyal
product.  By the above theorem, this notion can be extended to an
arbitrary formal deformation of $(\R^{2n},\omega)$ via the isomorphism
from statement (2), and this notion is well-defined by statement (3).
Finally, using Darboux's theorem, we can define the canonical trace on
a formal deformation of an arbitrary symplectic manifold
$(M^{2n},\omega)$ by requiring that the pull-back of such a functional
$\mathrm{tr}:\fmh\rightarrow\hbar^{-n}\Ch$ with respect to a
symplectic embedding of an open subset of $\R^{2n}$ into $M$ is a
canonical trace on (an open subset of) $\R^{2n}$.
\begin{prop}[\cite{NT}]
  For every formal deformation $(\fmh,\cdot)$ of a symplectic manifold
  $(M,\omega)$ the canonical trace exists and is unique.
\end{prop}

We denote the canonical trace by $\trc$.  Now we can formulate the
Fedosov-Nest-Tsygan index theorem \cite{Fe,NT} for a compact symplectic
manifold $M$. It relates the two objects defined above:
\begin{equation}
  \label{eq:FNT}
 \trc(1) = \int_M e^\theta\ahat(M).  
\end{equation}

This is an analog of the Grothendieck-Riemann-Roch formula, which have
already used in the previous section. For a compact symplectic
manifold $(M,\omega)$ with integral symplectic form, let $L$ be a line
bundle whose first Chern class is $\omega$. If $M$ has a compatible
K\"ahler structure and $L$ is endowed with an appropriate holomorphic
structure, then, for large $k$, the GRR theorem gives the following
expression for the dimension of the space of sections of $L^k$:
\begin{equation}
\label{eq:GRR}
\dim H^0(L^k) = \int_M e^{kc_1(L)}\todd(M). 
\end{equation}
This expression is an integer valued polynomial which only depends on
the symplectic structure. Note that the RHS, which in the projective
algebraic case is also known as the Hilbert polynomial, can be defined
for any symplectic manifold.

An alternative expression of the same type can be given using the
$\hat A$-genus:
\begin{equation*}
 \dim H^0(L^k)  = \int_M e^{kc_1(L)+c_1(M)/2}\ahat(M). 
\end{equation*}
Define the shifted Hilbert polynomial as
\begin{equation}
\label{eq:ahat}
P_L(k) = \int_M e^{kc_1(L)}\ahat(M). 
\end{equation}
Now define a deformation $(M,\cdot)$ of a symplectic manifold {\em
  basic} if $\theta(M,\cdot) = \omega/\hbar$. Then combining the above
formulas we obtain that given a basic deformation of a compact
symplectic manifold $(M^{2n},\omega)$ and a line bundle $L$ on $M$
with $c_1(L)=\omega$, we have \cite{Fe,NT}
\begin{equation}
\label{eq:bla}
  \trc(1) = P_L(\hbar^{-1}).
\end{equation}

\subsection{Algebraic manifolds and local deformations}
\label{sec:alg-loc}

For the purposes of our paper we need to reformulate the theory we
just outlined.

Let $A_0$ be a subalgebra of $\cim$ which separates points and is not
in the kernel of any non-zero complex differential operator on $M$.
Let $\Ah$ be an $\hbar$-adically complete associative algebra over
$\Ch$, such that $\Ah/\hbar\Ah\isom A_0$. Assume that there exists a
section $s:A_0\rightarrow \Ah$ of the natural map $\Ah\rightarrow A_0$
which, when extended to $A_0\hh$, gives an isomorphism
$s_\hbar:A_0\hh\simeq \Ah$. Then the formula $f\cdot_s g =
s^{-1}(s(f)s(g))$ defines a product on $A_0\hh$. If this product is
local (see the definition at the start of this section), then we say
that $s$ is a \emph{local} section.
\begin{deff}
  An algebra $\Ah$ which \emph{has} a local section is called a
  \emph{local deformation} of $A_0$.
\end{deff}

As it is clear from the construction above, the data $(\Ah,s)$ of an
algebra with a local section such that $A_0=\cim$ is equivalent to
that of a formal deformation of the manifold $M$.
\begin{lemma}
  The action of the gauge group defined in \eqref{eq:gauge}
  corresponds to the action $s\rightarrow s_\gamma =
  s\circ\gamma^{-1}$.
\end{lemma}
The proof is clear.  This means that the classification of formal
deformations of symplectic manifolds up to gauge transformations is
equivalent to the classification of local deformations of the algebra
$\cim$. Then the following statement is immediate:
\begin{prop}
  The Poisson structure, the characteristic class and the canonical
  trace $\Ah\rightarrow\C\hh$ induced by a pair $(\Ah,s)$ is
  independent of $s$.
\end{prop}
The interest in this statement lies in the possibility of studying
these invariants in examples of local deformations $A_\hbar$, which do
not have natural sections.

Such examples could arise in the following setting. Suppose that
$(M,\omega)$ is a smooth {\em real} affine algebraic manifold with a
symplectic Poisson structure on the space of algebraic functions.
Denote by $A_0$ the space of \emph{complex} algebraic functions on
$M$. While this is somewhat arbitrary, we choose our deformation ring
to be 
\begin{equation}
  \label{eq:dq}
  \dq = \{\text{Rational functions in }\, q,\text{ with no poles at }
  q=1\text{ and }0<|q|<1\}
\end{equation}

A non-commutative algebra $A_q$ over $\dq$, which is given by finite
number of generators and relations, such that $A_q/(q-1)A_q\isom A_0$
is often called a \emph{$q$-deformation} of $A_0$. It is a \emph{local
  $q$-deformation} if, in addition, the $\hbar$-adic completion
$A_\hbar=A_q\hat\tensor_{\Cq}\Ch$, where $q=e^{\ul\hbar}$, is a local
deformation of $A_0$ in the sense defined above. Assuming that the
Poisson structure derived from $\Ah$ is the given symplectic one, we
can ask the following questions:
\begin{enumerate}
\item What is the characteristic class of $A_q$, which is defined as
  the characteristic class of $\Ah$?
\item Can the canonical trace be defined on $A_q$?
\end{enumerate}

The first question is clear, although such computations are very
difficult. The second question requires some comment. One could hope
to start with a cyclic functional $\trq:A_q\rightarrow\dq$ and by
Taylor expansion obtain a functional $\tr_\hbar:\Ah\rightarrow\Ch$,
and then after choosing a local section $s$, one could arrive at a
functional $\tr_\hbar:A_0\hh\rightarrow\Ch$.  Next, one needs to show
that such a functional extends as a cyclic functional to $\cim\hh$ and
then finally, one could see if $\hbar^{-n}\tr_\hbar$ is the canonical
trace or not. Of course, it could very well happen that $A_q$ does not
have any cyclic functionals at all.

\textsc{An Example.} The only treatable, but non-trivial example with
which we are familiar is the quantum torus. The 2-dimensional torus
can be written written as an affine variety with generators
$\{U^\pm,V^\pm\}$ and relations $\{U^+U^-=V^+V^-=1\}$. It has a
canonical translation invariant symplectic form $\omega$, which we
normalize so that $\int\omega=1$. The non-commutative deformation has
the same, but now non-commutative, generators and relations and an
additional relation $UV=qVU$. Clearly, the monomials $\{U^mV^n|\,
m,n\in\Z\}$ form a basis of this algebra and this gives the necessary
section.
\begin{prop}
  1. This deformation is local, and the derived Poisson structure is
   the one, corresponding to $\omega$.\\
  2. The functional $\trq(U^mV^n)=\delta_m\delta_n$ is cyclic.  \\
  3. When passing to $\C\hh$, $q=e^{\ul\hbar}$, this functional is
  given by integration against $\omega$ and it is $\hbar$ times the
  canonical trace. \\
  4. The characteristic class of this deformation $\omega/\hbar$, thus
  the deformation is basic.
\end{prop}
We leave the proof of the proposition as an exercise to the reader.
The first two statements are easy, and the third one is doable. Note
that statement (3) and the F-NT index theorem imply statement (4).

One may generalize this setup, by starting with a cyclic functional on
$A_q$ with values in $\mathrm{Mer}^{\mathrm{asym}}_1(q)$, the
meromorphic functions on the unit disc $\{|q|<1\}$, which have an
asymptotic expansion in $q-1$ as $q$ approaches $1$ along the real
axis. The function $T_g(x;q)$ introduced in \eqref{eq:tqsut} is an
example of such a function. Then one can pass to values in $\C\hh$ by
setting $q=e^{\ul\hbar}$ and taking the asymptotic expansion
$\ase:\mathrm{Mer}^{\mathrm{asym}}_1(q)\rightarrow\C\hh$. As it happens, we
will need such a generalization.

In this paper, we construct $q$-deformations of the moduli spaces of
flat connections $\M^G[\Sigma_P](\mathbold{\sigma})$, and study the
questions raised above for these non-commutative algebras. We discuss
our results in detail in \secref{sec:conclusion}.

\section{G-Colored Graphs}
\label{sec:graph}

\subsection{Functions on the moduli space}
\label{sec:def:graph}

In this section we introduce \emph{$G$-colored graphs} which are a way
to represent the set $F(\M^G)$ of algebraic functions on $\M^G$.  They
are based on the notion of {\em graph connections} of Fock and Rosly
\cite{FR}, which was introduced as a discretization of gauge theory on
Riemann surfaces.  Our version is substantially similar but
technically more flexible than the original one. A similar
construction can be found in \cite{AMR}, where ordinary circular
holonomies are used and an almost identical notion appears in
\cite{PR}. Still there are some technical differences and the notion
of ``equivalence'' seems to be new. Our proofs will be brief since
they are similar to the original ones.

We keep the notation of \secref{sec:modspaces}.

A {\em $G$-colored } graph $f$ on $\SMP=\Sigma\setminus P$ consists of
\begin{itemize}
\item an oriented, not necessarily connected graph $\Gamma(f)$
  immersed in $\SMP$. We denote the set of edges of $\Gamma(f)$ by
  $\egam f$, the set of its vertices by $\vgam f$;
\item a coloring of each edge $e\in\egam f$ by a representation
  $C_f(e)$ of the group $G$.
\item a coloring of each vertex $v\in \vgam f$ by an invariant
  \[\phi_f(v)\in
  \left(\underset{e\rightarrow v}\tensor C_f(e)^*
    \underset{e\leftarrow v}\tensor C_f(e)\right)^G
\]
where the tensor products are taken over the incoming and outgoing
edges correspondingly.
\end{itemize}

\begin{rem}
  Note that the definition of coloring of the vertices above is
  somewhat imprecise, since we did not specify the order in which the
  tensor products are taken. Naturally, we are taking advantage of the
  fact that the space of invariants of a tensor product of
  representations of $G$ are naturally isomorphic to the space of
  invariants of the same representations tensored in a different
  order. In fact, using the orientation of the surface, we have a
  natural {\em cyclic orientation} of the edges adjacent to a vertex,
  but this still does not provide us with an ordering. This might seem
  like hair splitting here, but we will have to return to this
  question in the next section.
\end{rem}

Given an immersed graph $\Gamma$ in $\SMP$ with colored edges
$e\mapsto C(e)$, and a connection $\nabla$ on the trivial $G$-bundle
over $\SMP$, one can construct an element $\nabla_\Gamma\in
\underset{e\in E_\Gamma}\tensor \Hom(C(e),C(e))$ by taking the parallel
transports of $\nabla$ along the edges of $\Gamma$. Then given a
$G$-colored graph $f$, we obtain a number $f(\nabla)$ by pairing
$\nabla_\Gamma$ with $\underset{v\in \vgam f}\tensor\phi_{f}(v).$

\begin{lemma}
  \begin{enumerate}
  \item The number $f(\nabla)$ does not change if $\nabla$ is replaced
    by a gauge-equivalent connection.
  \item Let $f$ and $g$ be $G$-colored graphs and define $f\cup g$ to
    be the $G$-colored graph with $\Gamma_{f\cup
      g}=\Gamma(f)\cup\Gamma(g)$ and with coloring inherited from $f$
    and $g$. Then $(f\cup g)(\nabla) = f(\nabla)g(\nabla)$.
  \item If $\nabla$ is flat then $f(\nabla)$ is invariant under
    homotopic changes of the immersion of $\Gamma(f)$.
  \end{enumerate}
\end{lemma}

The proofs are straightforward and will be omitted.

Denote by $F(G,\Sigma_{P})$ the free vector space generated by all
$G$-colored graphs modulo homotopy of the embedding of the graph. Then
taking unions of the underlying graphs as in the Lemma endows
this space with an algebra structure. According to the statement (1) of
the Lemma, every element of this space defines a function on
$\M^G[\SMP]$. These functions are clearly algebraic and we have a
homomorphism $F(G,\Sigma_P)\rightarrow F(\M^G[\SMP])$ to the space of
algebraic functions on $\M^G[\SMP]$.

To describe the kernel of this map, fix an open disc $D\subset\SMP$
with boundary $\db$ a smooth embedded circle. Consider a $G$-colored
graph $f$ in generic position with respect to $D$.  This means that
every edge can intersect $\db$ only transversally and only finitely
many times. Then we can define a new, contracted $G$-colored graph
$f|D$, where the underlying immersed graph the quotient graph
$\Gamma/D\cap\Gamma$. This has a single vertex in $D$ and the edges
adjacent to this vertex correspond to the points of intersection of
the edges of $\Gamma$ with $\db$. The colorings of the edges and
vertices outside $D$ are inherited from $\Gamma$, while the coloring
of the new vertex can be obtained by contracting each edge $e$ in $D$
using the canonical diagonal element $\delta(C(e))$ in
$C(e)^{\star}\otimes C(e)$.  We will call two G-colored graphs,
related by such a contraction or a sequence of contractions, {\em
  equivalent}. Note that since we are free to move the graphs
homotopicly, the position of $D$ does not play any role.

Below we list a few important special cases of equivalent $G$-colored
graphs.
\begin{example}
\label{ex:first}
\begin{enumerate}
\item We can erase any edge which is colored by the trivial
  representation.
\item We can place a virtual 2-valent vertex colored by the diagonal
  element $\delta(V)$ on any edge colored by a representation $V$.
\item For two crossing edges, we can place a vertex at the
  intersection, colored by the permutation $P_{VW}:V\tensor
  W\rightarrow W\tensor V$.
\item A contractible loop colored by the representation $V_\lambda$ is 
  equivalent to the number $\dim V_\lambda$.
\item Fix an embedded interval $I\subset \Sigma_P$ and consider an
  $f\in F(G,\Sigma_P)$ in generic position with respect to $I$.
  For simplicity, assume that the $m$ edges intersecting $I$ are
  similarly oriented with respect to $I$; denote their colorings by
  $V_1,\dots,V_m$. Then $f$ is equivalent to $f\dagger I$ where
  $\Gamma(f\dagger I)$ is obtained from $\Gamma(f)$ by introducing two
  new $m+1$-valent vertices $S$ and $E$ joined by an edge $e$, as
  shown in Figure \ref{fig:cutting}. (The interval $I$ is represented
  by a thick line.) The element $f\dagger I\in
  F(G,\Sigma_P)$ is a sum of colored graphs which retain the colorings
  of those edges and vertices which are common with $f$, and thus have
  the form $g(f,V,\phi(S),\phi(E))$. Here $V=C(e)$, $\phi(S)\in
  \Hom(V_1\tensor\dots\tensor V_m,V)$ and $\phi(E)\in
  \Hom(V,V_1\tensor\dots\tensor V_m)$. These two spaces of invariants
  are naturally paired to each other via the formula
  $\tr_V(\phi(S)\phi(E))$, thus we can define
  \begin{equation}
    \label{eq:interval}
   f\dagger I =  \sum_{V\in\rg}\sum_i  g(f,V,\psi^i,\psi_i),
  \end{equation}
  where $\sum_i \psi^i\tensor\psi_i=\delta_\tr(\Hom(V_1\tensor\dots\tensor
  V_m,V))$ is the diagonal element induced by this pairing.
\end{enumerate}
\end{example}
\begin{rem}
  Note that the pairing $\tr(\phi(S)\phi(E))$ is somewhat
  redundant, since by Schur's Lemma we have 
  \begin{equation}
    \label{eq:schur}
\phi(S)\phi(E) = \frac{\tr_V(\phi(S)\phi(E))}{\dim V}\id_V.   
  \end{equation}
\end{rem}
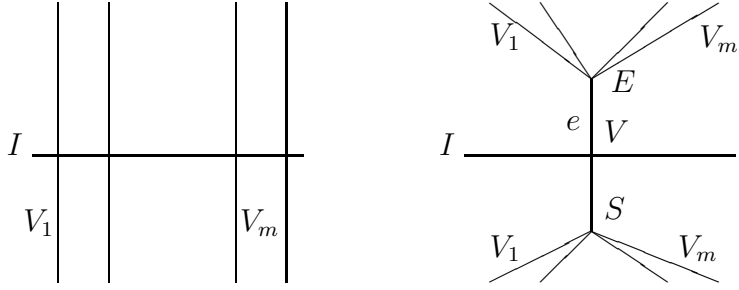
\begin{figure}
\setlength{\unitlength}{2100sp}%
\begingroup\makeatletter\ifx\SetFigFont\undefined%
\gdef\SetFigFont#1#2#3#4#5{%
  \reset@font\fontsize{#1}{#2pt}%
  \fontfamily{#3}\fontseries{#4}\fontshape{#5}%
  \selectfont}%
\fi\endgroup%
\begin{picture}(8712,3324)(901,-4273)
\put(8026,-2011){\makebox(0,0)[lb]{\smash{\SetFigFont{12}{14.4}{\rmdefault}{\mddefault}{\updefault}
$E$
}}}
\put(7951,-3511){\makebox(0,0)[lb]{\smash{\SetFigFont{12}{14.4}{\rmdefault}{\mddefault}{\updefault}
$S$
}}}
\put(7951,-2611){\makebox(0,0)[lb]{\smash{\SetFigFont{12}{14.4}{\rmdefault}{\mddefault}{\updefault}
$V$
}}}
\thinlines
\put(7201,-961){\line( 2,-3){600}}
\put(7801,-1861){\line( 1, 1){900}}
\put(6601,-4261){\line( 2, 1){1200}}
\put(7801,-3661){\line( 5,-2){1500}}
\put(7801,-3661){\line(-1,-1){600}}
\put(7801,-3661){\line( 3,-2){900}}
\put(7801,-1861){\line( 0,-1){1800}}
\put(1501,-961){\line( 0,-1){3300}}
\put(2101,-961){\line( 0,-1){3300}}
\put(4201,-961){\line( 0,-1){3300}}
\put(3601,-961){\line( 0,-1){3300}}
\put(4201,-961){\line( 0,-1){3300}}
\put(1501,-961){\line( 0,-1){3300}}
\put(1501,-961){\line( 0,-1){3300}}
\put(1096,-3661){\makebox(0,0)[lb]{\smash{\SetFigFont{12}{14.4}{\rmdefault}{\mddefault}{\updefault}$V_1$}}}
\put(3651,-3661){\makebox(0,0)[lb]{\smash{\SetFigFont{12}{14.4}{\rmdefault}{\mddefault}{\updefault}$V_m$}}}
\put(6601,-3961){\makebox(0,0)[lb]{\smash{\SetFigFont{12}{14.4}{\rmdefault}{\mddefault}{\updefault}$V_1$}}}
\put(8826,-3961){\makebox(0,0)[lb]{\smash{\SetFigFont{12}{14.4}{\rmdefault}{\mddefault}{\updefault}$V_m$}}}
\put(6601,-1461){\makebox(0,0)[lb]{\smash{\SetFigFont{12}{14.4}{\rmdefault}{\mddefault}{\updefault}$V_1$}}}
\put(9076,-1486){\makebox(0,0)[lb]{\smash{\SetFigFont{12}{14.4}{\rmdefault}{\mddefault}{\updefault}$V_m$}}}
\put(7501,-2461){\makebox(0,0)[lb]{\smash{\SetFigFont{12}{14.4}{\rmdefault}{\mddefault}{\updefault}$e$}}}
\put(6001,-2761){\makebox(0,0)[lb]{\smash{\SetFigFont{12}{14.4}{\rmdefault}{\mddefault}{\updefault}$I$}}}
\put(901,-2761){\makebox(0,0)[lb]{\smash{\SetFigFont{12}{14.4}{\rmdefault}{\mddefault}{\updefault}$I$}}}
\put(6601,-961){\line( 4,-3){1200}}
\put(7801,-1861){\line( 5, 3){1500}}
\put(2101,-961){\line( 0,-1){3300}}
\thicklines
\put(6301,-2761){\line( 1,0){3200}}
\put(1201,-2761){\line( 1, 0){3200}}
\end{picture}

\caption{\label{fig:cutting}$f$ and $f\dagger I$}
\end{figure}
Denote by $\mathrm{Col}(\Gamma)$ the set of all possible colorings of
the edges of a fixed immersed graph $\Gamma$ by \emph{irreducible}
representations, and by $\Phi(\Gamma,C)$ the linear space of colorings
of the vertices of a graph $\Gamma$ with a fixed coloring
$C\in\mathrm{Col}(\Gamma)$ of its edges. According to the following
proposition the kernel of the map discussed above is generated by our
notion of equivalence.
\begin{prop}
  \label{thm:fmgp}
  \begin{enumerate}
  \item Two equivalent $G$-colored graphs $f$ and $f|D$ take the same
    values on any flat connection $\nabla$.
  \item Let $\Gamma\subset\Sigma_P$ be an {\em embedded} graph, such
    that each {\em face} of $\Gamma$ is contractible and contains
    exactly one puncture. Then
\[ F(\M^G[\SMP]) = \underset{C\in \mathrm{Col}(\Gamma)}\oplus
\Phi(\Gamma,C).  \]
\item The kernel of the map $F(G,\Sigma_P)\rightarrow F(\M^G[\SMP])$
  is linearly generated by equivalence.
\end{enumerate}
\end{prop}
\begin{deff}  Define the graphs satisfying the condition in the 2nd
  statement {\em exact}. Such graphs always exist as long as there is
  at least one puncture. For an exact graph $\Gamma\in\Sigma_P$ define
  the \emph{dual graph} $\check\Gamma$ to be a graph embedded into
  $\Sigma$ with vertices at the punctures (the set $P$), and faces
  containing exactly one vertex of $\Gamma$ each.
\end{deff}
\emph{Proof}: (1). If we trivialize $\nabla$ over $D$, then the
parallel transport of $\nabla$ along the edges will be all equal to
the identity element of $G$. Then performing the partial contractions
in the definition of $f(\nabla)$ by contracting along the edges in $D$
only,
we arrive at $f|D(\nabla)$.\\
(2). This was pointed out in \cite{FR}. The statement follows from the
Peter-Weyl Theorem, since $\M^G$ is simply a product of groups divided
by the diagonal adjoint action. Note that coloring each edge by the
trivial representation and the vertices by the trivial invariants
gives the
unit element of the algebra.\\
(3). It follows from (2) that it is sufficient to prove that given an
exact graph $\Gamma$ any element $f\in F(G,\SMP)$ is equivalent to a
sum of $G$-colored graphs with underlying graph $\Gamma$.  Let
$\check\Gamma$ be the dual graph as described in the definition above.
Put the given  $G$-colored graph $f$  into general position with
respect to $\check\Gamma$. Then by performing the $\dagger$ operation
on $f$ with respect to each edge of $\check\Gamma$ we obtain a new
$G$-colored graph $f_{\check\Gamma}$ (or a sum of such) which
intersect each edge of $\check\Gamma$ exactly once and which is
equivalent to $f$.  Finally, using equivalence again, we can replace
each of the vertices in each face of $\check\Gamma$ by a single
vertex, thus obtaining a sum of $G$-colored graphs, each having
exactly one vertex in each face of $\check\Gamma$.  $\Box$
  
For each puncture $p\in P$ and representation $V$ define an element
$c_V^p\in F(G,\Sigma_P)$ with $\Gamma(c_V^p)$ a small counterclockwise
oriented circle around $p$ colored by $V$. We will also use the
notation $c_\lambda^p$ when $V=V_\lambda$. The graph underlying the
product $c_V^p c_W^p,$ is the union of two small concentric circles
around $p$.  The colored graph $(c_V^p c_W^p)\dagger I$, where $I$ is
an interval which intersects the two circles transversally, has 4
vertices. Contracting the two pairs of vertices joined by two edges,
it is easy to see that $c_V^p c_W^p$ is equivalent to $c_{V\otimes
  W}^p$. Thus the correspondence $c_V^p\mapsto V\in R(G)$ extends to a
homomorphism (in fact, isomorphism) of algebras.

\subsection{The Poisson structure}
\label{sec:Poisson}

Now we define a Poisson structure on the space of $G$-colored graphs
\cite{FR, AMR, PR}.  Fix an element $t\in(\sym^2(\lig))^G$. Then
for two $G$-colored graphs $f$ and $g$ in general position, and a
point $m\in \Gamma(f)\cap\Gamma(g)$ we can define a new $G$-colored
graph $f\cup^t_m g$, which is obtained from $f\cup g$ by placing a
vertex at $m$, colored by $P_{12}\circ t:C(e_f(m))\tensor C(e_g(m))
\rightarrow C(e_g(m))\tensor C(e_f(m))$, where $e_f(m)$ and $e_g(m)$
are the two edges containing $m$ and $P_{12}$ is the permutation
operator.

Now define
\begin{equation}
  \label{eq:def:poisson}
  \po fg = \sum_{m\in \Gamma(f)\cap\Gamma(g)}\sign(e_f(m),e_g(m))
  f\cup^t_m g,
\end{equation}
where the sign is obtained by comparing the orientations of the
ordered pair $(e_f(m),e_g(m))$ to the orientation of $\Sigma$.

\begin{prop}
\label{thm:poisson-center}
  \begin{enumerate}
  \item The operation $\po fg$ is well-defined on $\fmgp$, i.e. it is
    compatible with homotopy and equivalence.
  \item The operation $\po fg$ is a Poisson bracket on $\fmgp$.
  \item The elements $c_V^p$, $p\in P$ generate a subalgebra in the
    Poisson center of $\fmgp$ in the sense that $\po{c_V^p}{f}=0$ for
    all $f\in\fmgp$.
  \item The spaces $\M^G(\vec t)\subset\M^G[\Sigma_P]$ are symplectic
    leaves of this Poisson structure, and the induced symplectic
    form is exactly $\omega_{\vec t}$ (cf. \secref{sec:modspaces}).
  \end{enumerate}
\end{prop}
{\em Proof}: (1). To prove compatibility with homotopy, it sufficient
to show the property shown on Figure \ref{fig:poissonhomotopy}, where
vertices colored by $t$ are marked by circles.
\begin{figure}
\setlength{\unitlength}{2447sp}%
\begingroup\makeatletter\ifx\SetFigFont\undefined%
\gdef\SetFigFont#1#2#3#4#5{%
  \reset@font\fontsize{#1}{#2pt}%
  \fontfamily{#3}\fontseries{#4}\fontshape{#5}%
  \selectfont}%
\fi\endgroup%
\begin{picture}(7287,2760)(1426,-2986)
\thinlines
\put(3301,-1186){\circle{150}}
\put(1651,-1861){\line( 4, 1){1270.588}}
\put(2926,-1561){\line( 0,-1){1200}}
\put(2926,-1561){\line( 1, 1){937.500}}
\put(2626,-511){\line( 1,-1){1387.500}}
\put(4051,-1861){\line( 0, 1){  0}}
\put(4051,-1861){\line( 0, 1){  0}}
\put(7100,-1711){\circle{150}}
\put(7801,-2386){\circle{150}}
\put(6526,-1861){\line( 4, 1){1270.588}}
\put(7801,-1561){\line( 0,-1){1200}}
\put(7801,-1561){\line( 1, 1){937.500}}
\put(6676,-1261){\line( 1,-1){1387.500}}
\put(8101,-2611){\line( 0, 1){  0}}
\put(8101,-2611){\line( 0, 1){  0}}
\put(5026,-1636){\makebox(0,0)[lb]{\smash{\SetFigFont{12}{14.4}{\rmdefault}{\mddefault}{\updefault}=}}}
\put(2476,-361){\makebox(0,0)[lb]{\smash{\SetFigFont{12}{14.4}{\rmdefault}{\mddefault}{\updefault}$3$}}}
\put(6526,-1186){\makebox(0,0)[lb]{\smash{\SetFigFont{12}{14.4}{\rmdefault}{\mddefault}{\updefault}$3$}}}
\put(1426,-1936){\makebox(0,0)[lb]{\smash{\SetFigFont{12}{14.4}{\rmdefault}{\mddefault}{\updefault}$1$}}}
\put(2851,-2986){\makebox(0,0)[lb]{\smash{\SetFigFont{12}{14.4}{\rmdefault}{\mddefault}{\updefault}$2$}}}
\put(6301,-1861){\makebox(0,0)[lb]{\smash{\SetFigFont{12}{14.4}{\rmdefault}{\mddefault}{\updefault}$1$}}}
\put(7726,-2986){\makebox(0,0)[lb]{\smash{\SetFigFont{12}{14.4}{\rmdefault}{\mddefault}{\updefault}$2$
}}}
\end{picture}

\caption{\label{fig:poissonhomotopy}}
\end{figure}
This easily
follows from the identity
\[ t_{13}+t_{23} = (\Delta_{0}\otimes \id)(t)\in U(\lig)^{\tensor 3}, \]
where $\Delta_{0}$ is the coproduct in the universal enveloping
algebra $U(\lig)$. The indices, as usual, mark the embedding of
$V^{\tensor 2}$ into $V^{\tensor 3}$, e.g. $t_{13} =
P_{23}(t\tensor\id)$.  The compatibility with equivalence follows from
this because using homotopy it can always be arranged that there are
no intersection points in $D$. \\
(2). This statement follows from simple combinatorics of the
intersection points  (cf. \cite{AMR,PR}). \\
(3). Again, using homotopy we can arrange that an arbitrary graph does
not intersect a small circle around $p$. \\
(4). This statement is one of the main results of \cite{FR}.$\Box$

\subsection{The Poisson trace}
\label{sec:P-trace}

\begin{figure}
\setlength{\unitlength}{2100sp}%
\begingroup\makeatletter\ifx\SetFigFont\undefined%
\gdef\SetFigFont#1#2#3#4#5{%
  \reset@font\fontsize{#1}{#2pt}%
  \fontfamily{#3}\fontseries{#4}\fontshape{#5}%
  \selectfont}%
\fi\endgroup%
\begin{picture}(8712,3324)(901,-4273)
\put(8026,-2011){\makebox(0,0)[lb]{\smash{\SetFigFont{12}{14.4}{\rmdefault}{\mddefault}{\updefault}
$E$
}}}
\put(7951,-3511){\makebox(0,0)[lb]{\smash{\SetFigFont{12}{14.4}{\rmdefault}{\mddefault}{\updefault}
$S$
}}}
\put(7951,-2611){\makebox(0,0)[lb]{\smash{\SetFigFont{12}{14.4}{\rmdefault}{\mddefault}{\updefault}
$V$
}}}
\thinlines
\put(7201,-961){\line( 2,-3){600}}
\put(7801,-1861){\line( 1, 1){900}}
\put(6601,-4261){\line( 2, 1){1200}}
\put(7801,-3661){\line( 5,-2){1500}}
\put(7801,-3661){\line(-1,-1){600}}
\put(7801,-3661){\line( 3,-2){900}}
\multiput(7801,-1861)(0.00000,-116.12903){16}{\line( 0,-1){ 58.065}}
\put(1501,-961){\line( 0,-1){3300}}
\put(2101,-961){\line( 0,-1){3300}}
\put(4201,-961){\line( 0,-1){3300}}
\put(3601,-961){\line( 0,-1){3300}}
\put(4201,-961){\line( 0,-1){3300}}
\put(1501,-961){\line( 0,-1){3300}}
\put(1501,-961){\line( 0,-1){3300}}
\put(1096,-3661){\makebox(0,0)[lb]{\smash{\SetFigFont{12}{14.4}{\rmdefault}{\mddefault}{\updefault}$V_1$}}}
\put(3651,-3661){\makebox(0,0)[lb]{\smash{\SetFigFont{12}{14.4}{\rmdefault}{\mddefault}{\updefault}$V_m$}}}
\put(6601,-3961){\makebox(0,0)[lb]{\smash{\SetFigFont{12}{14.4}{\rmdefault}{\mddefault}{\updefault}$V_1$}}}
\put(8826,-3961){\makebox(0,0)[lb]{\smash{\SetFigFont{12}{14.4}{\rmdefault}{\mddefault}{\updefault}$V_m$}}}
\put(6601,-1461){\makebox(0,0)[lb]{\smash{\SetFigFont{12}{14.4}{\rmdefault}{\mddefault}{\updefault}$V_1$}}}
\put(9076,-1486){\makebox(0,0)[lb]{\smash{\SetFigFont{12}{14.4}{\rmdefault}{\mddefault}{\updefault}$V_m$}}}
\put(7501,-2461){\makebox(0,0)[lb]{\smash{\SetFigFont{12}{14.4}{\rmdefault}{\mddefault}{\updefault}$e$}}}
\put(6001,-2761){\makebox(0,0)[lb]{\smash{\SetFigFont{12}{14.4}{\rmdefault}{\mddefault}{\updefault}$I$}}}
\put(901,-2761){\makebox(0,0)[lb]{\smash{\SetFigFont{12}{14.4}{\rmdefault}{\mddefault}{\updefault}$I$}}}
\put(6601,-961){\line( 4,-3){1200}}
\put(7801,-1861){\line( 5, 3){1500}}
\put(2101,-961){\line( 0,-1){3300}}
\thicklines
\put(6301,-2761){\line( 1,0){3200}}
\put(1201,-2761){\line( 1, 0){3200}}
\end{picture}

\caption{\label{fig:cutting0}Cutting along $I$}
\end{figure}
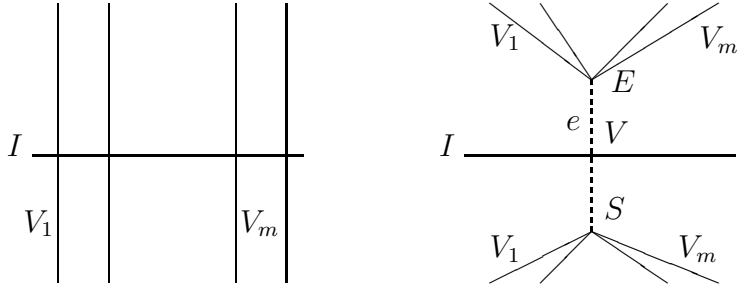

\propref{thm:fmgp} (2) allows us to define an augmentation
$H_\Gamma:F(\M^G[\SMP] )\rightarrow\C$ by projecting onto the exact
graph $\Gamma$ colored by trivial representations and invariants.  We
also give a more constructive formula for $H_\Gamma(f)$ for a general
$f\in F(G,\SMP)$ which will be useful for computations later on.
First, we define a variant of the operation $\dagger I$, denoted by
$\dagger_0 I$ and called {\em cutting}, which is similar to $\dagger$
with the difference that  $V$ is allowed to be the trivial representation
only.  Thus \eqref{eq:interval} is modified by
\begin{equation}
\label{eq:interval0}
 f\dagger_0 I =   \sum_i  g(f,\C,\psi^i,\psi_i),
\end{equation}
where $\sum_i
\psi^i\tensor\psi_i=\delta_\tr(\Hom(V_1\tensor\dots\tensor V_m,\C))$.
\begin{rem}
  We present a schematic picture of the cutting operation on Figure
  \ref{fig:cutting0}.  Note that now we can erase the edge between $S$
  and $E$, since it is colored by the trivial representation. However,
  as a mnemonic for the diagonal element $\delta_\tr$ that is
  inserted, we join the two vertices by a dashed line (chord).
\end{rem}
Now assume that $\Gamma$ is exact, and suppose that the graph
$\Gamma(f)$ is in a generic position with respect to the dual graph
$\check\Gamma$. By applying the cutting operation with respect to each
edge of $\check\Gamma$ we obtain a new $G$-colored graph
\[f\dagger_0\check\Gamma = 
f\prod_{e\in E_{\check\Gamma(f)}}\dagger_0 e.\] The graph underlying
$f\dagger_0\check\Gamma$ is a union of disjoint pieces, each located
on some contractible face.  Thus $f\dagger_0\check\Gamma$ is
equivalent to a number.
\begin{prop}
\label{thm:H-class}
  \begin{enumerate}
  \item The operation $f\dagger_0\check\Gamma$ is well defined on
    $F(\M^G[\SMP])$, i.e. it is compatible with homotopy and
    equivalence.
  \item The $G$-colored graph $f\dagger_0\check\Gamma$ is equivalent
    to the number $H_\Gamma(f)$.
  \item $H_\Gamma(f)$ does not depend on the choice of the exact graph
    $\Gamma$.
  \item $H_\Gamma(f)$ is given by a smooth top form on the smooth part
    of $\M^G[\Sigma_P]$.
  \item For any $f,g$, $H(\{f,g\})=0$.
  \end{enumerate}
\end{prop}
{\em Proof}. (1). Clearly compatibility with equivalence holds with
respect to any disc $D$ which lies entirely in one of the faces of
$\check\Gamma$. Using this we may assume that each face contains only
one vertex. Then the only relevant homotopy relation is moving one of
these vertices across one of the cutting edges. We leave proving this
case as an exercise to the reader. \\
(2). Because of part (1) we may assume that $\Gamma(f)=\Gamma$, in
which case the
statement is obvious. \\
(3). It is easy to see that any two exact graphs are related by the
operation of contracting single edge in some exact graph $\Gamma$.
This corresponds to removing an edge in $\check\Gamma$, and the
statement now follows from compatibility with equivalence. Hence
from now on we can omit the index $\Gamma$ in $H_\Gamma$.\\
(4). Note that there is another simple way to define $H_\Gamma$. From
our earlier discussion it is clear that the space $\M^G$ is a quotient
of a product of copies of the group $G$, corresponding to the edges of
$\Gamma$ by the action of a product of copies of the group $G$
corresponding to the vertices of $\Gamma$. Since the only matrix
coefficient on a compact Lie group whose integral does not vanish is
that of the trivial representation, we see that the measure induced on
$\M^G$ by the operation $H_\Gamma$ is simply the push-forward of the
Haar measure on the product of the groups. This push-forward is
clearly
smooth whenever the action is locally trivial.\\
(5). The proof of this statement is analogous to that of
\propref{thm:quantum-trace} (cf. \remref{thm:poisson-commute}). It
follows from the fact that a self-intersecting edge colored by
$V_\lambda$, with the tensor $P_{12}\circ t$ inserted at the
intersection point is equivalent to an edge without self-intersection,
colored the same way and multiplied by $-2C(\lambda)$, where
$C(\lambda)$ is the value of the Casimir operator, normalized using
$t$, on $V_\lambda$. This, in turn, follows from the equality
$t=\Delta C-1\tensor C-C\tensor 1$. There is a somewhat exotic proof
of this statement in \cite{PR}. $\Box$

Consider now the following general situation.  Assume that
$\pi:M\rightarrow N$ is a fibration between
two compact, smooth manifolds, and assume that the manifolds are
endowed with smooth volume forms $\mu_M$ and $\mu_N$ such that the
volume of both manifolds is 1. Then one can define the push-forward
operation on continuous functions $\pi_*:C^0(M)\rightarrow C^0(N)$ by
integrating with respect to the natural measure $\mu_M/\pi^*\mu_N$
along the fibers. The proof of the following formulas will be omitted:
\begin{lemma}
  Let $\{c_i,c^i\}_{i=0}^\infty$ be dual bases of functions
  on $N$, i.e.  $\int_N c^ic_j\,\mu_N=\delta_{ij}$ and the functions
  $\{c_i\}$ are complete in $L^2(N,\mu_N)$. Then
  \begin{enumerate}
  \item 
\begin{equation}
    \label{eq:pushforward}
    \pi_*(f) = \sum_i\left(\int_M f\pi^{*}(c^i)\,\mu_M\right)\,c_i
  \end{equation}
\item  The permanence equation
\begin{equation}
  \label{eq:permanence}
  \pi_*(f\pi^*(g)) = \pi_*(f)g
\end{equation}
holds for any $f\in C^0(M)$ and $g\in C^0(N)$.
  \end{enumerate}
\end{lemma}

For notational simplicity we will concentrate on the $|P|=1$ case from
here on.  Now set $N=\orb_{\mathrm{reg}}(G)$ with measure induced by
the Haar measure on $G$, $\pi=\holp$ and $M\subset\M^G[\Sigma_p]$ the
smooth part of $\holp^{-1}(N)$ with the measure defined by $H$ above.
While the technical conditions of the Lemma do not hold, the
conclusions do (cf. \cite[Theorem 4.2]{BL}), thus we can conclude
\begin{prop}
\label{thm:rat-series}
  Define a functional $\tr$ on $\fmgp$ by the series
\[ \tr(f) = \sum_{\lambda\in \Omega^+} H(fc^\lambda)\bar{\chi}_\lambda. \]
\begin{enumerate}
\item Then the series converges pointwise at every regular orbit and
  for every $G$-colored graph $f$, takes values in continuous
  functions on $\orb_\mathrm{reg}(G)$.
\item  The $\C$-valued functional
$\tr(f)^\sigma$ obtained by evaluation at a non-special conjugacy
class $\sigma$ can be obtained by integration along $\M^G(\sigma)$
with respect to a smooth measure $\mu_\sigma$.
\item We have
 \[ \tr(fc_\lambda) = \tr(f)\chi_\lambda. \]
\end{enumerate}
\end{prop}

Note that for a non-special $t$, we  have a symplectic form
$\omega_t$ on $\fmgp$ which also produces a smooth volume form. The
following equality is due to  Witten \cite{Wi,KL,BL} 
\[ c(g,G)\delta(t)\mu_t =
\frac{\omega(t)^{\dim\M^G(\sigma)}}{\dim\M^G(\sigma)!} \] 
where $c(g,G)$ is defined by \eqref{eq:w-def}. 

To demonstrate the power of the ``cutting calculus'' described above,
we finish this section with the computation of $\tr(1)$ for an
arbitrary number of punctures. As
\[ \tr(f) = \sum_{\vec\lambda}H\left(f\prod_{p\in P}
  c^p_{\vec\lambda(p)}\right)\prod_{p\in P}\bar\chi_{\vec\lambda(p)}\]
computing $\tr(1)$ involves computing $H(\prod_{p\in
  P}c^p_\lambda(p))$. Choose an exact graph $\Gamma$ and cut the graph
under consideration by the edges of $\check\Gamma$. Then every edge
$e$ of $\check\Gamma$ will cut through two edges with opposite
orientations.  In order for them to give a non-zero contribution, they
have to have the same coloring. Thus $H(\prod_{p\in P}c_\lambda(p))=0$
unless $\lambda(p)=\lambda(p')=\lambda$ for $p,p'\in P$. Then
according to \eqref{eq:schur} the contribution at each cutting edge is
a factor of $(\dim V_\lambda)^{-1}$, while the remaining graph
consists of a union of loops colored by $\lambda$, one on each face of
$\check\Gamma$.  These faces correspond to the vertices of $\Gamma$
and each contributes a factor of $\dim V_\lambda$ (cf. Example
\ref{ex:first} (2)). Thus we obtain
\[ H\left(\prod_{p\in P}c^p_\lambda\right) = (\dim
V_\lambda)^{|V_\Gamma|-|E_\Gamma|} = (\dim V_\lambda)^{2-2g-|P|}.\]
Here we used that the Euler characteristic of $\Sigma$ is $2-2g$ and
that the faces of $\Sigma$ are in one-to-one correspondence with the
punctures. Thus we have
\begin{equation}
  \label{eq:class-volseries}
  \tr^{\vec\sigma}(1) = \sum_{\lambda} \frac{\prod_{p\in
    P}\chi_{\lambda}(\vec\sigma(p))}{(\dim V_\lambda)^{2g-2+|P|}}.
\end{equation}

\section{Ribbon graphs and the moduli algebra}
\label{sec:ribbon}

In this section we construct a non-commutative $q$-deformation of the
algebra $\fmgp$, based on the representation theory of quantum groups.
This algebra is similar or equivalent to the ``moduli algebra''
constructed in \cite{AGS}and \cite{BR} (cf. also \cite{AMRq,Tu}). Our
construction is much more geometric and transparent, however, and the
calculations are much simpler. We clarify the relation of this algebra
to the ribbon categories of Reshetikhin and Turaev \cite{RT}, which
simplifies the contruction a great deal.

Note that since we will need to pass to special values of $q$, instead
of the standard quantum group $\ue_q(\lig)$, we will work over a
smaller, ``non-restricted'' algebra defined over the
ring $\dq$, which we introduced in \secref{sec:formdef} (cf.
\secref{sec:appendix} for some details). We will use the generic
symbol $\ue$ for this algebra.

We start with a technical prelude. In defining the quantum analog of
the $G$-colored graphs, we need a somewhat more geometric version of
the coloring of the vertices than that in \cite{RT}. The notion of
\emph{cyclic invariant} that we introduce allows us to define an
algebra of the ``correct'' size.

\subsection{The Reshetikhin-Turaev map} 
Recall the construction of ribbon categories of Reshetikhin and Turaev
\cite{RT}.  

Define a {\em band} as a rectangle embedded in oriented 3-space which
has its sides marked as follows: the starting edge, the ending edge,
the left edge and the right edge.  Alternatively, one can think of an
arrow drawn on one side (the ``marked side'') of the rectangle
parallel to one pair of edges.  In particular, the band and its
boundary segments are oriented. By attaching  an edge $e$ of a band
to an oriented segment $I$, we mean that restricting the embedding of the
band to $e$ is an orientation-reversing embedding of $e$ into $I$.

Let $\ue$ be an appropriate non-restricted Ribbon Hopf algebra algebra
\cite{RT,CP} (cf. \secref{sec:q-basic} and remark above).  Fix
$\lstart$ and $\lend$, two parallel oriented lines in $\R^3$. A {\em
  ribbon configuration} is a union of bands of two types, ribbons and
coupons, which projects into the strip between the two lines and such
that the lines, ribbons and coupons are all disjoint except that the
starting and ending edges of each ribbon are attached to either a line
or to the starting or ending edge of a coupon. If we associate an
irreducible representation of $\ue$ to each ribbon of a particular
ribbon configuration $\mbC$, then the whole configuration, as well as
each coupon $\mbc$ acquires a {\em type}: $(S(\mbc),E(\mbc))$ which is
simply the list of the colorings of the ribbons at the starting and
the ending edges of the coupon, with the direction of the arrows
recorded as follows: we record a $+$ if the arrow points towards the
coupon and a $-$ if it points away from the coupon. Thus
$(S(\mbc),E(\mbc))$ has the form
$([(V_1,\eps_1),\dots,(V_k,\eps_k)],[(W_1,\rho_1),\dots,(W_l,\rho_l)])$,
where $\eps_i=\pm 1, \rho_i=\pm1$ and $V_i,W_j\in\rue$.

A \emph{$\ue$-ribbon configuration} is one where each ribbon is
colored by a representation of $\ue$, and each coupon by a
$\ue$-invariant map
\begin{equation}
    \label{eq:uinv}
   V_s(\mbc) = V_1^{\eps_1}\tensor\dots \tensor V_k^{\eps_k} \rightarrow
 W_1^{\rho_1}\tensor\dots \tensor W_l^{\rho_l}=V_e(\mbc)
\end{equation}
with the convention that if the coloring is $V$ then on the starting
end $V^+=V$ and $V^-=V^*$ where the latter is the left representation
of $\ue$ defined on the dual space to $V$ using the antipode; on the
other end it is the other way around: $W^-=W$ and $W^+=W^*$.

A fundamental result of \cite{RT} is that each $\ue$-ribbon
configuration defines a morphism of $\ue$-representations
$\vstart(\mbC)\rightarrow\vend(\mbC)$, which depends on the ribbon
configuration up to homotopy only. Denote by $I_{(S,E)}$ the linear
space of intertwining maps \eqref{eq:uinv}, which we will call
invariants of type $(S,E)$.  Then the result may be summarized by
saying that a ribbon configuration with colored edges defines a
homotopy invariant product
\begin{equation}
  \label{eq:coupons}
\mathrm{RT}:  \tensor_{\mbc}I_{(S(\mbc),E(\mbc))}\rightarrow
I_{(S(\mbC),E(\mbC))}. 
\end{equation}
We will call this product the \emph{Reshetikhin-Turaev} map.

\subsection{Cyclic invariants and colored ribbon graphs}
\label{sec:cyclic}

The improvement that we suggest is the following: define a
\emph{cyclic type} $(T)^\cyc$ to be a cyclicly ordered set of
representations of $\ue$, each marked with a $+$ or a $-$.  Clearly,
each ordinary type induces a cyclic type by assigning a cyclic order
to the $V$'s and $W$'s the obvious way, listing the representations
counterclockwise, ignoring which ones are attached to the lower edge
and which ones to the upper edge.  If $(T)^\cyc$ has $m$
representations, there will be $m(m+1)$ different types corresponding
to $(T)^\cyc$.  Indeed, one can reduce a cyclic ordering on $m$
elements to an ordinary ordering in $m$ different ways, and then one
can divide an ordered set of $m$ elements into two ordered sequences
in $m+1$ ways. Denote the cyclic type derived from an ordinary type
$(S,E)$ by $(S,E)^\cyc$.

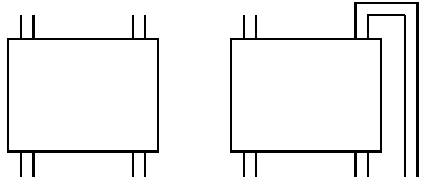
\begin{figure}
\setlength{\unitlength}{2047sp}%
\begingroup\makeatletter\ifx\SetFigFont\undefined%
\gdef\SetFigFont#1#2#3#4#5{%
  \reset@font\fontsize{#1}{#2pt}%
  \fontfamily{#3}\fontseries{#4}\fontshape{#5}%
  \selectfont}%
\fi\endgroup%
\begin{picture}(4974,2124)(1489,-2173)
\thinlines
\put(4201,-1861){\framebox(1800,1350){}}
\put(4351,-211){\line( 0,-1){300}}
\put(4501,-211){\line( 0,-1){300}}
\put(5701,-211){\line( 0,-1){300}}
\put(5851,-211){\line( 0,-1){300}}
\put(4351,-1861){\line( 0,-1){300}}
\put(4501,-1861){\line( 0,-1){ 75}}
\put(4501,-1936){\line( 0,-1){ 75}}
\put(4501,-2011){\line( 0,-1){ 75}}
\put(4501,-2086){\line( 0,-1){ 75}}
\put(4201,-1861){\framebox(1800,1350){}}
\put(5701,-1861){\line( 0,-1){ 75}}
\put(5701,-1936){\line( 0,-1){ 75}}
\put(5701,-2011){\line( 0,-1){ 75}}
\put(5701,-2086){\line( 0,-1){ 75}}
\put(4201,-1861){\framebox(1800,1350){}}
\put(5851,-1861){\line( 0,-1){ 75}}
\put(5851,-1936){\line( 0,-1){ 75}}
\put(5851,-2011){\line( 0,-1){ 75}}
\put(5851,-2086){\line( 0,-1){ 75}}
\put(1501,-1861){\framebox(1800,1350){}}
\put(1651,-211){\line( 0,-1){300}}
\put(1801,-211){\line( 0,-1){300}}
\put(3001,-211){\line( 0,-1){300}}
\put(3151,-211){\line( 0,-1){300}}
\put(1651,-1861){\line( 0,-1){300}}
\put(1801,-1861){\line( 0,-1){ 75}}
\put(1801,-1936){\line( 0,-1){ 75}}
\put(1801,-2011){\line( 0,-1){ 75}}
\put(1801,-2086){\line( 0,-1){ 75}}
\put(1501,-1861){\framebox(1800,1350){}}
\put(3001,-1861){\line( 0,-1){ 75}}
\put(3001,-1936){\line( 0,-1){ 75}}
\put(3001,-2011){\line( 0,-1){ 75}}
\put(3001,-2086){\line( 0,-1){ 75}}
\put(1501,-1861){\framebox(1800,1350){}}
\put(3151,-1861){\line( 0,-1){ 75}}
\put(3151,-1936){\line( 0,-1){ 75}}
\put(3151,-2011){\line( 0,-1){ 75}}
\put(3151,-2086){\line( 0,-1){ 75}}
\put(5701,-211){\line( 0, 1){150}}
\put(5851,-211){\line( 1, 0){450}}
\put(6301,-211){\line( 0,-1){1950}}
\put(5701,-61){\line( 1, 0){750}}
\put(6451,-61){\line( 0,-1){2100}}
\end{picture}

\caption{Operation $A_{er}$\label{fig:bendingaband}}
\end{figure}

There are operations of bending an edge (or ``attaching a candy
cane'') which map spaces of invariants of the same cyclic type into
each other. Given a type
\[(S, E)=
([(V_1,\eps_1),\dots,(V_k,\eps_k)],[(W_1,\rho_1),\dots,(W_l,\rho_l)]),
\]
and an intertwiner $\phi\in I_{(S,E)}$, let 
\[ A_{er}\phi=(\mathrm{id}\otimes a(W_l,\rho_l)) (\phi\otimes
\mathrm{id}_{ W_l^{-\rho_l}}),\]
where $a(W, \rho):W^{\rho}\otimes
W^{-\rho}\rightarrow \C$ is the standard coinvariant (cf. Appendix). 
Then $A_{er}\phi$ is an
invariant of type
$$
([(V_1,\eps_1),\dots,(V_k,\eps_k),
(W_l,\rho_l)],[(W_1,\rho_1),\dots,(W_{l-1},\rho_{l-1})])
$$
A pictorial representation of this operation is shown on Figure
\ref{fig:bendingaband}.  We can define three other operations
$A_{el},A_{sr}$ amd $A_{sl}$ which bend the edges $W_1,V_k$ and $V_1$
respectively. Clearly, have $A_{er}A_{sr}=\mathrm{id}$ and
$A_{el}A_{sl}=\mathrm{id}$.  Iterating these maps we get various maps
between these spaces of invariants.
\begin{lemma}
  The spaces of invariants of the same fixed cyclic type
\[ \{I_{(S,E)}|\,(S,E)^\cyc=(T)^\cyc\}\] 
are canonically isomorphic under these maps.
\end{lemma}
\emph{Proof}: This is implicit in the original paper \cite{RT}. One
needs to check a generalization of the relation $(\mathrm{Rel}_{13})$
of \cite{RT}. An instance of such a relation is 
\[ \phi(\mu v_1,\mu v_2,\dots,\mu v_l) = \phi( v_1, v_2,\dots, v_l) \]
for an invariant $\phi\in\mathrm{Hom}(V_1\tensor
V_2\tensor\dots\tensor V_l,\C)$.  This follows from the fact that
$\mu$ is a group-like element, and $\epsilon(\mu)=1$.$\Box$

We will call invariants related by the above maps \emph{cyclically
  equivalent}. The Lemma allows us to define the space $I_{(T)^\cyc}$
of \emph{cyclic $\ue$-invariants} of type $(T)^\cyc$ as equivalence
classes of intertwiners of the same cyclic type related by the above
maps.

\begin{example}
\label{triv-vert}
The simplest example of a cyclic invariant is the trivial invariant of
type $([(V,+)],[(V,-)])^\cyc$. This can be interpreted as the 6
different invariant maps, which correspond to the diagram (Fig.1) of
\cite{RT}, or equivalently to the 6 maps of page 167 of \cite{CP}
denoted $\iota_{V}^+, \iota_{V}^{-}, \alpha_V^{+}, \alpha_V^{-},
\beta_V^{+}, \beta_V^{-}.$
\end{example}

\begin{rem}
  \label{thm:cycred}
Since $\mu$ equals to the identity when $q=1$, cyclic equivalence of
quantum invariants corresponds to ordinary cyclic permutations and
contractions of the factors in the classical case.
\end{rem}

Now we can define a non-commutative generalization of $G$-colored
graphs.
\begin{deff}
  A {\em $\ue$-colored ribbon graph} $\fb$ consists of
\begin{itemize}
\item a {\em ribbon graph} $R(\fb)$, which is a union of oriented
  discs (vertices) and ribbons embedded into $\R^3$, The discs and
  ribbons are disjoint except that each end of each ribbon is attached
  to the boundary of a disc, as always, respecting the orientations. We
  denote the set of ribbons of $R(\fb)$ by $E_{R(\fb)}$ and the set of
  vertices by $V_{R(\fb)}$.
\item A coloring $C_\fb:E_{R(\fb)}\rightarrow\rue$ of the ribbons by
 representations of $\ue$. Clearly such a coloring
  assigns a cyclic type $(T(v))^\cyc$ to each vertex $v\in V_R(f)$
\item A coloring of each vertex $v\in V_R$ by a cyclic invariant
  $\phi(v)$ of type $(T(v))^\cyc$.
\end{itemize}
\end{deff}

This definition does not allow for free edges; we will stipulate their
existence whenever necessary. The key point is that the
Reshetikhin-Turaev map given in \eqref{eq:coupons} is compatible with
cyclic equivalence.
\begin{lemma}
  \label{lem:cycprod}
  Let $\fb$ be a $\ue$-colored ribbon graph with $m$ free edges
  embedded in the interior of a cylindrical surface $CS$ (e.g. in
  $\{x^2+y^2<1\}$) in $\R^3$ in such a way that the free edges are
  attached to a fixed circle (e.g. $\{x^2+y^2=1, z=0\}$) on $CS$. A
  choice of type, compatible with the cyclic type, induced by the
  colorings and the cyclic ordering of the $m$ free edges, as well as
  a choice of a type for each vertex, give rise to a $\ue$-ribbon
  configuration. The RT map \eqref{eq:coupons} induced by this
  configuration applied to the colorings of the vertices give rise to an
  invariant map between the tensor products of the colorings of the
  free edges. Then all the invariant maps thus obtained are cyclicly
  equivalent.
\end{lemma}
The proof of this lemma is left to the reader as an exercise.  It will
be important for us that this operation of replacing a ribbon graph by
a single vertex is defined over $\dq$. 

We can pair the free edges against a  vertex and obtain the
following:
\begin{cor}
\label{cor:rt-inv}
The algorithm described in the Lemma associates a well-defined number
to every $\ue$-colored ribbon graph (without free edges).
\end{cor}
This is our version of the Reshetikhin-Turaev invariant. A special
case is
\begin{cor}
  The pairing between $\mathrm{Hom}(\C,V_1\otimes \cdots \otimes V_l)$ and
  $\mathrm{Hom}(V_1\otimes \cdots \otimes V_l,\C)$ given by
  composition is cyclically invariant.
\end{cor}
One could also say that there is a canonical pairing
\[ \langle,\rangle:I_{(T)^\cyc} \tensor I_{(T^*)^\cyc} \rightarrow
\C \] between cyclic invariants if type $(T)^\cyc$ and those of the
dual type $(T^*)^\cyc$, which changes the cyclic orientation to the
opposite and changes each representation to its dual.

\subsection{Quantization of moduli spaces}
\label{sec:qmod}

The notions of cyclic invariants and colored ribbon graphs permit us
to $q$-deform the constructions of the previous section. The
constructions below are completely parallel to the classical case of
the previous section.

Let $\Sigma$ be a compact Riemann surface embedded into 3-space,
oriented outward, with a set of marked points $P$, as
before.  Denote by $\Se$ a small open neighborhood of $\Sigma$ and by
$\pi$ a projection $\pi:\Se\rightarrow\Sigma$ and let
$\smpe=\pi^{-1}(\Sigma\less P)$.  Then define $\fqgsp$ to be the free
$\dq$-module, generated by the homotopy classes of embeddings of
colored ribbon graphs into $\smpe\subset\R^3$. We can define a
{\em product} of two colored ribbon graphs $\fb,\gb\in\fqgsp$ as
the disjoint union of the two ribbon graphs $\fb'$ and $\gb'$, where
$\fb'$ is in the interior of $\Sigma$ and is homotopic to $\fb$, while
$\gb'$ is in the exterior of $\Sigma$ and is homotopic to $\gb$.  Note
that this is an associative but, generally, non-commutative product.
An analogous operation was used earlier by Turaev in \cite{Tu} (also
cf. \cite{AMRq}).

Now we define the notion of {\em equivalence}.  We fix a disc
$D\subset \SMP$ and a $\ue$-colored ribbon graph $\fb$ in $\smpe$ such
that $R(\fb)\cap \pi^{-1}(\db)\subset\Sigma$. This means that the
ribbons intersecting $\pi^{-1}(\db)$ are effectively attached to
$\db$. This is a version of the notion of \emph{generic position} with
respect to $D$. Then using \lemref{lem:cycprod} we can construct a new
graph $\fb|D\in\fqgsp$ by replacing the part $\pi^{-1} D\cap R(\fb)$
by a single vertex which can be arranged to lie entirely in $\Sigma$.
The colorings of the edges and vertices outside $D$ are inherited from
$\fb$, while the coloring of the new vertex is obtained
\lemref{lem:cycprod}.  Now define the algebra $F^q(\M^G[\Sigma_P])$ to
be the quotient of $\fqgsp$ by the $\dq$ linear subspace generated by
this notion of equivalence, which subspace is clearly also an ideal.
Thus $F^q(\M^G[\Sigma_P])$ is endowed with a natural associative
algebra structure over $\dq.$

Again we have the basic examples of equivalence:
\begin{example}
\label{thm:q-examps}
\begin{enumerate}
\item We can erase an edge which is colored by the trivial
  representation.
\item We can divide every ribbon into two pieces joined by a vertex
  colored by the trivial cyclic invariant (cf. Example
  \ref{triv-vert}).
\item A small ring-like ribbon contractible in $\Sigma^\epsilon_P$ and
  colored by $V_\lambda$ is equivalent to the $q$-number $q\dim
  V_\lambda$.
\item We can replace two overcrossing (resp. undercrossing) ribbons
  colored by $V$ and $W$, by 4 ribbons attached to a coupon colored by
  $P_{VW}R^{(+)}_{VW}$ (resp. $P_{VW}R^{(-)}_{VW}$), where as usual
  $R^{(+)}$ denotes the $R$ matrix and $R^{(-)}$ denotes
  $R_{21}^{-1}.$ (Figure \ref{fig:overcrossing}).
\begin{figure}
\setlength{\unitlength}{2947sp}%
\begingroup\makeatletter\ifx\SetFigFont\undefined%
\gdef\SetFigFont#1#2#3#4#5{%
  \reset@font\fontsize{#1}{#2pt}%
  \fontfamily{#3}\fontseries{#4}\fontshape{#5}%
  \selectfont}%
\fi\endgroup%
\begin{picture}(4074,1674)(1039,-1723)
\thinlines
\put(1351,-1711){\line( 3, 4){450}}
\put(1051,-1711){\line( 1, 0){300}}
\put(1051,-1711){\line( 3, 4){612}}
\put(2251,-61){\line(-3,-4){450}}
\put(2251,-61){\line( 1, 0){300}}
\put(2551,-61){\line(-3,-4){612}}
\put(1051,-61){\line( 3,-4){1224}}
\put(2251,-1711){\line( 1, 0){300}}
\put(2551,-1711){\line(-3, 4){1224}}
\put(1351,-61){\line(-1, 0){300}}
\put(1051,-61){\line( 3,-4){1224}}
\put(2251,-1711){\line( 1, 0){300}}
\put(2551,-1711){\line(-3, 4){1224}}
\put(1351,-61){\line(-1, 0){300}}
\put(3601,-61){\line( 1, 0){300}}
\put(3901,-61){\line( 3,-4){225}}
\put(4801,-61){\line(-3,-4){225}}
\put(3751,-1411){\framebox(1200,1050){}}
\put(3826,-361){\line(-3, 4){225}}
\put(3601,-61){\line( 1, 0){300}}
\put(3901,-61){\line( 3,-4){225}}
\put(4576,-361){\line( 3, 4){225}}
\put(4801,-61){\line( 1, 0){300}}
\put(5101,-61){\line(-3,-4){225}}
\put(3826,-1411){\line(-3,-4){225}}
\put(3601,-1711){\line( 1, 0){300}}
\put(3901,-1711){\line( 3, 4){225}}
\put(4576,-1411){\line( 3,-4){225}}
\put(4801,-1711){\line( 1, 0){300}}
\put(5101,-1711){\line(-3, 4){225}}
\put(4176,-886){\makebox(0,0)[lb]{\smash{\SetFigFont{12}{14.4}{\rmdefault}{\mddefault}{\updefault}
$PR^{(+)}$
}}}
\put(2851,-961){\makebox(0,0)[lb]{\smash{\SetFigFont{12}{14.4}{\rmdefault}{\mddefault}{\updefault}
$\sim$
}}}
\end{picture}
\caption{\label{fig:overcrossing}}
\end{figure}\\
\item The operation $\dagger I$, as well as the operation of {\em
    cutting} $\dagger_0 I$ are defined as in the previous section. The
  pairing $\tr_V(\phi(S)\phi(E))$ needs to be replaced by
  $\tr_V(\mu\phi(S)\phi(E))$ as it is natural in the theory of quantum
  groups (cf. \secref{sec:q-basic}).  Just as in the definition of
  equivalence, the condition of generic position is $\pi^{-1}(I)\cap
  R(\fb)\in I$.
\item Assume that a ribbon $\mathbf{e}$ of a $\ue$-colored ribbon
  graph $\fb$ is colored by a representation $V_\lambda$. Then we have
  $\tilde\fb=v(\lambda)\fb$, where $\tilde\fb$ is $\fb$ with the
  ribbon $\mathbf{e}$ twisted according to the orientation of 3-space
  by 360 degrees, and $v(\lambda)$ is the value of the (central)
  ribbon element $v$ in the representation $V_{\lambda}$ (cf.
  \cite{RT}, \secref{sec:q-basic}).  Twisting in the opposite
  direction induces multiplication by $v(\lambda)^{-1}$.
\end{enumerate}
\end{example}

For an embedded graph $\Gamma\subset\Sigma_P$, there is a well-defined
ribbon graph $R_\Gamma\subset\Sigma_P$ obtained by thickening $\Gamma$
in $\Sigma$. In particular, one can associate canonical elements
$\bc^p_\lambda$ to the Poisson central elements $c^p_\lambda$ defined in the
previous section.

According to \remref{thm:cycred}, there is a well-defined classical
limit of a cyclic invariant to a classical invariant, extending
$\ev:\dq\rightarrow\C$.
\begin{prop}
\label{thm:fqb}
  \begin{enumerate}
  \item Define a map over $\C$
 \[   \mathrm{red}_F:F^q(\M^G[\Sigma_P])\rightarrow F(\M^G[\Sigma_P]),
\]
by shrinking the width of the ribbons to 0, projecting them onto
$\Sigma_P$ and applying the abovementioned reduction of the cyclic
invariants to the classical ones. This map is an algebra homomorphism
over $\C$.
  \item The space of all colorings of an exact graph $\Gamma$ span
    $F^q(\M^G[\Sigma_P])$.
  \item The elements $\{\bc^p_\lambda|\,p\in P,\,\lambda\in
    \Omega^+\}$ are in the center of $F^q(\M^G[\Sigma_P])$. They span an 
    algebra isomorphic to $R(G)^{\oplus|P|}\tensor\dq$.
  \end{enumerate}
\end{prop}
\emph{Proof}: (1). This amounts to checking that the relations in the
quantum case reduce to the classical ones when
$q=1$. This is straightforward.\\
(2.) The proof is similar to that of \propref{thm:fmgp} (2). By
iterating the $\dagger e$ operation with respect to the edges to the
dual graph $\check\Gamma$, one can write any element
of $\fqmgp$ as a sum of elements with underlying graph $\Gamma$.\\
(3.) The statement is clear since one can move a small circle around
the line $\pi^{-1}(p)$ for some puncture $p$, the graph underlying
$\bc^p_V$, past any other ribbon graph using a homotopy. Moreover,
using the same argument as in the classical case, in $\fqmgp$ we have
$\bc^p_{V\otimes W}=\bc^p_{V}\bc^p_{W}=\bc^p_{W}\bc^p_{V}$.$\Box$

The definition of the analog of the operation $H_\Gamma$ is a bit more
subtle, because out of the 3 definitions we gave in the commutative
case (projection, cutting, integration) only the cutting operation is
clearly well-defined.  
\begin{prop}
\label{thm:quantum-trace}
  \begin{enumerate}
  \item The functional $H_\Gamma$ defined by the cutting operation,
    does not depend on $\Gamma$. Thus we have a well-defined
    functional $H:F^q(\M^G[\Sigma_P])\rightarrow \dq$. Also, the
    functionals $H$ in the classical and quantum cases are compatible
    with the reduction map $\mathrm{red}_F$, i.e. we have
    $H(\mathrm{red}_F)(\fb))=\ev(H(\fb))$.
  \item The Poisson structure induced on $F(\M^G)$ by the evaluation map
    $\mathrm{red}_F$ and the relation $q=e^{\pi i\hbar}$ coincides
    with the one defined in \secref{sec:Poisson}.
  \item For $\fb,\gb\in F^q(\M^G[\Sigma_P])$, we have
    $H(\fb\gb)=H(\gb\fb)$
 \end{enumerate}
\end{prop}
\emph{Proof}: (1). The proof is the same as in the classical case.\\
(2). This can be derived from the results \cite{FR}. We will the
details here.\\
(3). It is sufficient to prove the statement for the case when
$R(\fb)$ is the thickening of an exact graph $\Gamma$ and $R(\gb)$ is
the same graph, but with edges oriented in the opposite way. Choose an
ordering of the edges at each vertex of $\Gamma$ compatible with the
cyclic order. Every edge of the dual graph $\check\Gamma$ is crossed
by two edges of $\fb\gb$ and $\gb\fb$ which have the same coloring but
opposite orientation. To perform the cutting operation, we move these
two edges side by side, so that we
achieve the condition
\[ \pi^{-1}(\check\Gamma)\cap R(\fb\gb) = \check\Gamma\cap
R(\fb\gb), \] and we do the same for $\gb\fb$. Note that this involves
making a non-canonical choice, but we make the same choice in both
cases. After performing the cutting with respect to the edges of
$\check\Gamma$ we obtain that
\[ H(\fb\gb) = \prod_{e\in E_\Gamma}(\qdim\, C(e))^{-1} \prod_{v\in V_\Gamma}
\langle \phi_{\fb}(v),\phi_{\gb}(v)\rangle_1, \] where
$\langle,\rangle_1$ is a certain pairing between the corresponding
invariants which depends on the particular choice that we made at each
edge of $\check\Gamma$. The formula for $H(\gb\fb)$ is the same, but
with a pairing $\langle,\rangle_2$ replacing $\langle,\rangle_1$. The
difference between the two pairings is that in the first case, the
ribbon graph $R(\gb)$ is above $R(\fb)$, and in the second case it is
below. The two cases are related by the homotopic move of rotating by
360 degrees the piece of $R(\gb)$ remaining on a particular face after
the cutting, so that it ends up under the corresponding piece of
$R(\fb)$. Thus the difference between $H(\fb\gb)$ and $H(\gb\fb)$ will
be a twist of $\pm 360$ degrees at every edge of $\check\Gamma$,
wherever the two pieces of graphs are joined. At an edge of
$\check\Gamma$, which we assume to be colored by $V_\lambda$, this
twist contributes a factor of $v^\pm(\lambda)$ on one side and
$v^\mp(\lambda)$ on the other, which cancel each other. This ends the
proof. $\Box$
\begin{rem}
  \label{thm:poisson-commute}
  Parts (2) and (3) together imply \propref{thm:H-class} (5), but the
  proof of (3) given above has a simple semiclassical version giving a
  proof of the Poisson trace property of $H$. Instead of the 360
  degree twists, in that case one encounters self-intersecting edges
  with the operator $P_{12}\circ t$ (cf. \secref{sec:Poisson}) inserted at
  the point of self-intersection. 
\end{rem}

\section{The trace}
\label{sec:proptraces}

Now we are ready to define the fixed holonomy quantized moduli spaces.
Again to avoid further complicating our notation, we will assume that
$G$ is simply-laced.

Fix a set of regular conjugacy classes
$\vec\sigma:P\rightarrow\orb_{\mathrm{reg}}(G)$ and define the
quotient by the ideal
\[ F^q(\M^G(\vec{\sigma})) = F^q(\M^G[\Sigma_P])/
\langle\{\bc_\lambda^p =
\chi_\lambda(\vec\sigma(p))|\, p\in P,\, \lambda\in \Omega^+\}\rangle 
\]

Define a functional $\trq$ on $\fqmgp$ by the series
\begin{equation}
\label{eq:trq}
 \trq^{\vec\sigma}(\fb) = \sum_{\vec\lambda}
H\left(\fb \prod_{p\in P}
\bc_{\vec\lambda(p)}^p\right)\prod_{p\in
P}\bar{\chi}_{\vec\lambda(p)}(\vec\sigma(p)). 
\end{equation}
Formally, the series takes values in functions in $q$ and $|P|$ copies 
of (a completion of) $R(G)$.

Our main result is formulated in the Theorem below. As we mentioned
earlier, we are only proving this statement for $G=\sut$ in this
paper, although several partial results are proved in the general case.

We will use the term ``convergence'' in the punctured unit disc, to
mean absolute and uniform convergence of holomorphic functions on each
ringlike domain $\{\epsilon\leq|q|\leq1-\epsilon\}$.

\begin{thm}
\label{thm:main}
Let $\fb\in\fqmgp$, $G=\sut$ and $\vec\sigma$ as above. 
\begin{enumerate}
\item Then the series \eqref{eq:trq} defining $\trq^{\vec\sigma}(\fb)$
converges to a holomorphic function on the punctured unit disc.
\item For every $p\in P$, $\lambda\in\Omega^+$
\[ \trq(\fb\bc^p_\lambda) = \trq(\fb)\chi_\lambda \]
thus the evaluation $\trq^{\vec\sigma}$ descends to the quotient
$F^q(\M^G(\vec{\sigma}))$.
\item If $\vec\sigma$ is not special (cf. \secref{sec:sev-pun}), then
  $\trq^{\vec\sigma}(\fb)$ has an asymptotic expansion at $q=1$. More
  precisely, there is a function $\tr_\hbar(\fb)$, analytic in a
  neighborhood of 0, and positive constants $\tau, C$ such that:
\[ \vert  \trq(\fb)-\tr_\hbar(\fb)\vert <C
e^{-\frac{\tau}{|\hbar|}}, \] where $q=e^{\pi i \hbar}$ and $\hbar\in
i\R^+$ is sufficiently small. Moreover, $\tr_\hbar(\fb)$ is a rational 
function in $q$ and $\hbar$.
\end{enumerate}
\end{thm}
\begin{rem}
  The notation $\tr_\hbar(\fb)$ is somewhat inconsistent. We will
  sometimes use $\ase(\trq(\fb))$ instead.
\end{rem}
The proofs of parts (1) and (3) are given after
\propref{thm:weyl-sym}. \\
\emph{Proof of }(2):  Consider the diagonal element $\sum_\lambda
\chi_\lambda\tensor\bar\chi_\lambda$ with respect to the standard
quadratic form on $R(G)$: $(\alpha,\beta)=H_G(\alpha\beta)$, where
$H_G$ is the projection onto the trivial character. This form
is manifestly invariant: $(\alpha\gamma,\beta)=(\alpha,\beta\gamma)$,
thus the diagonal element has a similar property:
\[\sum_\lambda \chi_\lambda\chi_\mu\tensor\bar\chi_\lambda=\sum_\lambda
\chi_\lambda\tensor\chi_\mu\bar\chi_\lambda.
\]
This implies the statement since the algebra spanned by
$\{\bc_\lambda\}$ is simply another copy of the algebra $R(G)$.$\Box$

\subsection{The one puncture case}
\label{sec:onepun}
Again, first we consider the $|P|=1$ case.  Let $p$ be a marked point
on a surface $\Sigma$ of genus $g$, and fix a conjugacy class $\sigma$
of the group $G$.

Our goal is to study the infinite series
\[ \trq(\fb) = \sum_{\lambda \in P^+} H(\fb \bc^p_{\lambda})
\chi_{\lambda}(\sigma).\]

We first express $H(\fb \bc^p_{\lambda})$ in terms of intertwiners of
irreducible finite dimensional representations of $\ue.$

Let $\Gamma$ be the usual exact graph with vertex $o$, $2g$ edges,
$\{e_i\}_{i=1}^{2g}$, and a single face containing the puncture $p$.
Then the ribbon graph $\Gamma(\bc^p_{\lambda})$ is the thickening of
a small circle around $p$, which is homotopic to the product
$e_1e_2e_1^{-1}e_2^{-1}\dots e_{2g-1}e_{2g}e_{2g-1}^{-1}e_{2g}^{-1}$
taken in $\pi_1(\Sigma_p)$.

By \propref{thm:fqb} (2), any colored ribbon graph $\fb\in\fqmgp$ is
equivalent to one with underlying ribbon graph $R_\Gamma$, the
thickening of the graph $\Gamma$, thus we can assume that, in fact,
$R(\fb)=R_\Gamma$.  Then $\fb$ is given by the colorings of the edges:
$C_{\fb}(e_{2i-1})=V_{\mu(i)}, C_{\fb}(e_{2i})=V_{\nu(i)}$,
$i=1,\dots,2g$, and an invariant
$\phi_{\fb}\in\mathrm{Hom}_{{\ue}}(\mathcal{V},\dq)$, where we
have denoted by $\mathcal{V}=\tensor_{i=g}^1 \mathcal{V}_i,$ with
$\mathcal{V}_i= V_{\nu(i)}^{*}\tensor V_{\mu(i)}^{*}\tensor
V_{\nu(i)}\tensor V_{\mu(i)} .$

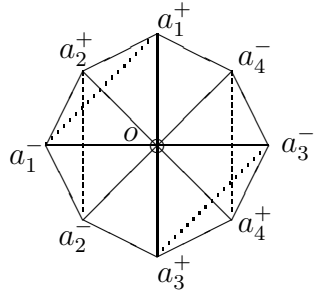
\begin{figure}
\setlength{\unitlength}{2047sp}%
\begingroup\makeatletter\ifx\SetFigFont\undefined%
\gdef\SetFigFont#1#2#3#4#5{%
  \reset@font\fontsize{#1}{#2pt}%
  \fontfamily{#3}\fontseries{#4}\fontshape{#5}%
  \selectfont}%
\fi\endgroup%
\begin{picture}(3225,3210)(1576,-2911)
\thinlines
\put(3301,-1261){\line(1,1){900}}
\put(3301,-1261){\line(-1,1){900}}
\put(3301,-1261){\line(1,-1){900}}
\put(3301,-1261){\line(-1,-1){900}}
\put(3301,-1261){\line(0,1){1350}}
\put(3301,-1261){\line(0,-1){1350}}
\put(3301,-1261){\line(1,0){1350}}
\put(3301,-1261){\line(-1,0){1350}}
\put(3301, 89){\line(-2,-1){900}}
\put(3301, 89){\line(2,-1){900}}
\put(1951,-1261){\line(1, 2){450}}
\put(1951,-1261){\line(1,-2){450}}
\put(2401,-2161){\line(2,-1){900}}
\put(3301,-2611){\line(2, 1){900}}
\put(4201,-2161){\line(1, 2){450}}
\put(4651,-1261){\line(-1, 2){450}}
\multiput(2401,-361)(0,-97){18}{\line(0,-1){54}}
\multiput(1951,-1261)(73,73){18}{\line(0,1){36}}
\multiput(4201,-361)(0,-97){18}{\line(0,-1){54}}
\multiput(3301,-2611)(73,73){18}{\line(0,1){36}}
\put(3301,164){\makebox(0,0)[lb]{\smash{\SetFigFont{12}{14.4}{\rmdefault}{\mddefault}{\updefault}$a_1^+$}}}
\put(3301,-1261){\circle{150}}
\put(2111,-200){\makebox(0,0)[lb]{\smash{\SetFigFont{12}{14.4}{\rmdefault}{\mddefault}{\updefault}$a_2^+$}}}
\put(1506,-1486){\makebox(0,0)[lb]{\smash{\SetFigFont{12}{14.4}{\rmdefault}{\mddefault}{\updefault}$a_1^-$}}}
\put(2106,-2461){\makebox(0,0)[lb]{\smash{\SetFigFont{12}{14.4}{\rmdefault}{\mddefault}{\updefault}$a_2^-$}}}
\put(3301,-2911){\makebox(0,0)[lb]{\smash{\SetFigFont{12}{14.4}{\rmdefault}{\mddefault}{\updefault}$a_3^+$}}}
\put(4276,-2311){\makebox(0,0)[lb]{\smash{\SetFigFont{12}{14.4}{\rmdefault}{\mddefault}{\updefault}$a_4^+$}}}
\put(4801,-1336){\makebox(0,0)[lb]{\smash{\SetFigFont{12}{14.4}{\rmdefault}{\mddefault}{\updefault}$a_3^-$}}}
\put(4276,-286){\makebox(0,0)[lb]{\smash{\SetFigFont{12}{14.4}{\rmdefault}{\mddefault}{\updefault}$a_4^-$}}}
\put(2888,-1210){\makebox(0,0)[lb]{\smash{\SetFigFont{12}{14.4}{\rmdefault}{\mddefault}{\updefault}$o$}}}
\end{picture}

\caption{\label{fig:cartwheel}Genus 2 case.}
\end{figure}

The dual graph $\Gc$ is isomorphic to $\Gamma$, has its vertex at $p$
and each of its $2g$ edges intersect exactly one edge of $\Gamma$ at
one point. Denote these points of intersection by
$\{a_i\}_{i=1}^{2g}$, correspondingly. We can represent $H(\fb
\bc^p_{\lambda})$ by performing the cutting operation with respect to
the edges of $\check\Gamma$. The resulting colored ribbon graph will
have the form of a cartwheel lying entirely above the face of
$\check\Gamma$.  If we represent the insertion of the diagonal element
by a dotted line, or \emph{chord}, between the two relevant vertices
of $H(\fb \bc^p_{\lambda})$, then we obtain the ``cartwheel with a
snow chain'' diagram with $4g+1$ vertices $\{o,a_i^\pm\}$ with $a_i^+$
and $a_i^-$ joined by a chord as shown schematicly on the Figure
\ref{fig:cartwheel} for the case $g=2$.

Next, for $\lambda,\mu\in P^+$ and a finite dimensional representation
$V$, introduce the notation
$$I(V;\lambda, \mu)=\Hom(V_\lambda, V_\mu\tensor V)\text{ and }
I^*(V;\lambda,\mu)=\Hom(V_\mu, V_\lambda\tensor V^*).$$
and
denote by $I(V;\lambda)$ the space $I(V;\lambda,\lambda).$ This
notation is appropriate since $I$ and $I^*$ are naturally dual to each
other. Indeed, let $\phi\in I(V;\lambda,\mu)$ and $\psi\in
I^*(V,\lambda,\mu)$. If we use the shorthand $\psi\circ\phi$ for the
composition $(\psi\tensor\id)\phi$ and $\langle\rangle_V$ for the natural
invariant pairing on $V^*\tensor V$, then the expression
$\langle\psi\circ\phi\rangle_V$ can be considered a constant, since it
represents an intertwiner from the irreducible representation
$V_\lambda$ to itself.  The resulting pairing is non-degenerate and
gives rise to the diagonal element $\delta(V;\lambda,\mu) \in
I^*(V;\lambda,\mu)\tensor I(V;\lambda,\mu)$ defined by
$\delta(V;\lambda,\mu)=\sum \alpha^i\tensor \alpha_i$ where
$\{\alpha_i\}$ and $\{\alpha^i\}$ are dual bases with respect to the
pairing, i.e.  $\langle\alpha_i\circ\alpha^j\rangle_V = \delta_{i}^{j}$.

If $V, W$ are two finite dimensional $\ue$-modules, let
\[\xi_{V,W}(\lambda)=\sum_{j,k} \alpha^{j}\circ \beta^{k}\circ
\alpha_{j}\circ \beta_{k}\in
\mathrm{Hom}(V_{\lambda},V_{\lambda}\otimes V^{*}\otimes W^{*}\otimes
V\otimes W),
\] 
where $\sum_{j} \alpha^{j}\otimes \alpha_{j}=\delta(V;\lambda)$ and
$\sum_{k} \alpha^{k}\otimes \alpha_{k}=\delta(W;\lambda)$. As usual,
we will write $\xi_{\lambda\mu}$ instead of $\xi_{V_\lambda V_\mu}$.
\begin{lemma}
\label{lem:cartwheel}
We have
\begin{equation}
\label{eq:cartwheel}
 H(\fb c^p_\lambda ) = \frac{1}{(\qdim
  V_\lambda)^{2g-1}}(\id_{V_{\lambda}}\otimes \phi_{\fb}) 
\xi_{\nu(g)\mu(g)}(\lambda)\circ\cdots \circ
\xi_{\nu(1)\mu(1)}(\lambda),
\end{equation}
where again $\circ$ stands for the composition of intertwiners as
above, each acting on $V_\lambda$ only.
\end{lemma}
\emph{Proof}: After cutting the cartwheel diagram between $a_{2g}^-$
and $a_1^+$ we obtain the expression on the RHS exactly, with each
chord contributing a factor $(\qdim V_\lambda)^{-1}$. When we make
this cut, we lose a factor of $\qdim V_\lambda$, hence the exponent
$1-2g$ in the expression.\fini

As a particular case of this expression, we can take $\fb=1$ and obtain:
\begin{equation}
H(c^p_{\lambda})=\frac{1}{(\qdim V_\lambda)^{2g-1}}
\end{equation}
for the $|P|=1$ case. This formula can be found in \cite{AGS,BR}. In
general, this leads to the quantum version of
\eqref{eq:class-volseries}:
\begin{equation}
  \label{eq:quant-volseries}
  \trq^{\vec\sigma}(1) = \sum_{\lambda} \frac{\prod_{p\in
    P}\chi_{\lambda}(\vec\sigma(p))}{(\qdim V_\lambda)^{2g-2+|P|}}.
\end{equation}
For $G=\sut$ this is exactly the series $\tilde T(\vex;q)$ introduced
in \secref{sec:qrat}, where $\vex$ and $\vec\sigma$ are related by the 
exponential map as usual.

Next, we would like to study the behavior of \eqref{eq:cartwheel} in
$\lambda$.  In order to track the $\lambda$-dependence of such
an expression, one would like to identify the spaces $I(V;\lambda)$ as
$\lambda$ varies. This can be done using the notions of ``expectation
value''and ``fusion matrices'' introduced in
\cite{ES,EV1}.  As we need a refinement of these objects, we give a
small introduction to the subject in \secref{sec:q-dynqg}, while only
providing the basic definitions below.

Our plan is to identify the spaces $I(V;\lambda,\mu)$ and the
corresponding intertwiner spaces for the Verma modules $\tilde
I(V;\lambda,\mu)$ with $V[\lambda-\mu]$ via the \emph{expectation value map}
$\langle\rangle$ (cf. \defref{thm:exp-val}).  This would
identify, for example, all of the spaces $I(V;\lambda)$ with the same
vector space $V[0]$, which would facilitate a universal treatment. To
this end we need the expectation value map to be an isomorphism, and
this turns out to be true most of the time. Indeed, if $\mu$ is
sufficiently far from the walls of the Weyl chamber relative to $\nu$,
i.e. $\mu$ is generic with respect to $\nu$ (cf.
\defref{thm:generic}), then the expectation value map is an
isomorphism between $I(V_{\nu};\lambda,\mu )$ and $V_{\nu}[\lambda-\mu
]$. In addition, in the Verma module case, the map $\langle\rangle$
has a right inverse $v\mapsto\phi^v_\lambda$.

For any sufficiently generic $\lambda$, this allows one to introduce
the \emph{fusion matrices} $J_{VW}(\lambda)\in\ndmf(V\tensor W)$,
defined by \eqref{eq:fusion-matrix}, which represent the operation of
composition of the intertwiners under the identification
$\langle\rangle$. This object has a natural generalization
$J_{UVW\dots}$ to tensor products with more than 2 factors. The operator
$Q_W(\lambda)\in\ndmf(W)$ is defined by taking the trace of
$J_{W*W}(\lambda)$ (cf. \eqref{eq:defq}). Finally, the \emph{dynamical
  Weyl group} elements $A_{V,w}(\lambda):V[\nu]\rightarrow V[w(\nu)]$,
defined for $\lambda\in\Omega^+$, $w\in W_G$, represent the standard
inclusion of Verma modules $M_{w.\lambda}\hookrightarrow M_\lambda$
under the same identification.

All objects defined in the previous paragraph have a universal form: 
$J(\lambda)=J_{12}(\lambda)$, $J_{1,2,\dots}$, $Q_(\lambda)$,
$A_w(\lambda)$ in an appropriate completions of tensor products of $\ue$.
The following lemma is a simple consequence of the notation we
introduced.
\begin{lemma}
  Let $V$ be a finite dimensional representation, $\lambda,\mu\in \Omega^+$
  and $\mu$ sufficiently generic. Then under the identification of
  $I(V;\lambda,\mu)$ with $V[\lambda-\mu]$ and $I^{*}(V;\lambda,\mu)$
  with $V[\lambda-\mu]^*$, the diagonal element
  $\delta(V;\lambda,\mu)$ corresponds to the element
\[ \sum_i
v^{i}\otimes Q_V(\lambda)^{-1}v_i\in V[\lambda-\mu]^*\tensor
V[\lambda-\mu],\] where $\sum_i v^i\tensor v_i =
\delta(V[\lambda-\mu])$.
\end{lemma}

The Lemma allows us to write down a formula for $H(\fb
c^{P}_{\lambda})$ in terms of fusion matrices only.
 \begin{prop}
\label{thm:Hfc=Rf}
   For notational convenience, denote by $W(i),i=1\cdots 4g$ the
   spaces
\[ V_{\nu(g)}^* [0], V_{\mu(g)}^*[0], V_{\nu(g)}[0],
V_{\mu(g)}[0],\cdots, V_{\nu(1)}^*[0], V_{\mu(1)}^*[0],
V_{\nu(1)}[0],V_{\mu(1)}[0]\] in that order. If $\lambda\in \Omega^+$ is
sufficiently generic, then we have:
\begin{equation}
\label{eq:Hfc=Rf}
 H({\bf f}\bc^{p}_{\lambda})=\frac{1}{(\qdim V_\lambda)^{2g-1}}
\phi_{\bf f}J_{1,\cdots,4g}(\lambda)\prod_{i=1}^{g}
Q_{W(4i-1)}^{-1}(\lambda)Q_{W(4i)}^{-1}(\lambda)\omega
\end{equation}
where $\omega=\otimes_{i=1}^{g} \delta(W(4i-1))_{4i-3,4i-1}\tensor
\delta(W(4i))_{4i-2,4i}$ is the diagonal element of
$\otimes_{i=1}^{4g} W(i).$
\end{prop}
The key point is that while the LHS of \eqref{eq:Hfc=Rf} is defined
for dominant integral weights $\lambda$, the RHS is meaningful for
arbitrary (generic) $\lambda$. This allows us to interpret the RHS in
representation theoretical terms for arbitrary $\lambda$.  
\begin{prop}
  \label{thm:q-dependence}
  Let $\fb\in
  F^q(\M^G)$. Then there exists a function $\hat R_{\fb}$, which is rational
  of non-positive degree in the variables $q_\alpha$,
  $\alpha\in\Delta^+$ with coefficients in $\dq$ and with poles along
  the ``hyperplanes'' $q_\alpha-q^m$, such that for
  $\lambda\in \Omega^+$, sufficiently generic, one has
  \begin{equation}
    \label{eq:vr}
   H(\fb \bc^p(V_\lambda))=\frac1{(\qdim
  V_{\lambda})^{2g-1}}\hat R_{\fb}\left(q_\alpha\mapsto
\frac{q^{(\alpha,\lambda)}-1}{q-1}, \alpha\in\Delta^+\right), 
  \end{equation}
where $\mapsto$ means substitution.
\end{prop}
\emph{Proof}: We only give a proof in the case of $G=\sut$. This
statement could be derived form results of \cite{EV1}, but in this
paper, since we are mostly working with $\sut$, we chose to give
explicit formulas from which the statement is manifest. We are hoping
that this will give a better idea of the complexity of the
computations to the reader. No such formulas are known in the case of
other groups.  

For $G=\sut$, the operators $J$ and $Q$ are given in a concrete basis
of each irreducible representation in \eqref{eq:Jofsu2} and
\eqref{eq:Qsu2}. Comparing these with \eqref{eq:Hfc=Rf} and the
formula \eqref{eq:higher-j} for $J_{1,\dots,N}$, the statement of the
proposition follows immediately.$\Box$

Denote the function obtained by substitution in the RHS of
\eqref{eq:vr} by $R_{\fb}(\lambda)$. Define the shifted Weyl group
action on $\tdual$ by $w.\lambda = w(\lambda+\rho)-\rho$ for all $w\in
W_G$.  We will now show that the function $R_{\fb}(\lambda)$ is
invariant under the shifted Weyl group action.  We will use the
properties of the operators $A_{w,V}$ defined above and in
\secref{sec:q-dynqg}. 
\begin{prop}
\label{thm:weyl-sym}
The function $R_{\fb}(\lambda)$ invariant under the shifted action of
the Weyl group: $R_{\fb}(\lambda)=R_{\fb}(w.\lambda)$ for every $w\in
W_G$. 
\end{prop}
\emph{Proof}: From the relation (\ref{eq:AJ=JA}), we obtain:
\begin{eqnarray}
&&\Delta^{4g}(A_w(\lambda))J_{1,\cdots,4g}(\lambda)=\hfill\\
&&J_{1,\cdots,4g}(w.\lambda)A_{w}^{(4g)}(\lambda)
A_{w}^{(4g-1)}(\lambda-h^{(4g)})\cdots
A_{w}^{(1)}(\lambda-h^{(2)}-\cdots h^{(4g)}),\nonumber 
\end{eqnarray}
where $\Delta^n$ is the iterated coproduct and $h^{(m)}$ means that
the element $h$ acts on the $m$th component of the tensor product.
Since $\phi_{\bf f}$ is an intertwiner, we have
$\phi_{\fb}\Delta^{(4g)}(A_{w}(\lambda))
=\epsilon(A_{w}(\lambda))\phi_{\fb}=\phi_{\fb}$.  Then we obtain:
\begin{equation}
R_{\fb}(\lambda)=\phi_{\fb}J_{1,\cdots,4g}(w.\lambda)
\prod_{j=4g}^{1}A_{w,W(j)}(\lambda)\prod_{i=1}^{g} 
Q_{W(4i-1)}^{-1}(\lambda)Q_{W(4i)}^{-1}(\lambda)\omega. 
\end{equation}
As a result it is sufficient to show the following relation:
$$A_w^{(4)}(\lambda)A_w^{(3)}(\lambda)A_w^{(2)}(\lambda)A_w^{(1)}
(\lambda)Q_{2}(\lambda)^{-1}Q_{1}(\lambda)^{-1}\omega=
Q_{2}^{-1}(w.\lambda)Q_{1}^{-1}(w.\lambda)\omega.$$
This is equivalent
to showing that
\begin{equation}
Q(w.\lambda)^{-1}\vert_{V[0]}=
A_{w}(\lambda)Q^{-1}(\lambda)S(A_{w}(\lambda))\vert_{V[0]}
\end{equation}
on the zero weight subspace.  This follows from the relation
(\ref{eq:Q=SAQA}).$\Box$ \smallskip
\begin{rem}
  Note that in the $\sut$ case $J_{VW}(\lambda)$ restricted to
  $V[0]\otimes W[0]$ has simples poles of the type
  $(1-q^{2(\lambda-k)})^{-1}$, $k\in \mathbb{N}.$ From
  (\ref{eq:QQ=v}), $Q_{V}(\lambda)^{-1}(\lambda)$ restricted to $V[0]$
  has simple poles of the type: $(1-q^{2(-\lambda-2-k)})^{-1}$, $k\in
  \mathbb{N}.$ This is a consistency check with the Weyl invariance of
  $R_{\fb}(\lambda)$ and \eqref{eq:Hfc=Rf}.
\end{rem}

We have shown that far from the walls of the dominant chamber $H(\fb
\bc_\lambda^p)$ has a certain form, which will be used to study
$\trq(\fb)$.  The expression for $\trq(\fb)$, however, contains a sum
over all dominant integral weights, and therefore we have to analyze
the behavior of $H(\fb \bc_\lambda^p)$ on the hyperplanes of non-generic
weights as well.  One needs to refine the results of \cite{EV1} in
order to be able to treat these cases, and we will not do this here.

Nevertheless, we can continue the study in the $G=\sut$ case where
this problem can be circumvented.  Indeed, as it was pointed out in
\remref{thm:r-generic} (3), in this case, there are only finitely many
non-generic weights for every representation $V$.

\emph{Proof of parts (1) and (3) of \thmref{thm:main}}.  Denote by $\nu$ the
fundamental weight of $\sut$. We now return to the notation 
of \secref{sec:3series}: assume that a number
$0<x<1/2$ represents the conjugacy class $\sigma$, so that
$\exp(x)=\xh\in\sigma$. We can conclude from
\propref{thm:q-dependence} and \propref{thm:weyl-sym} that
\[ \delta(\xh)\trq^\sigma(\fb) -\fel\sum_{n\in\Z}\!\!\!\!* 
\frac{e_{n\nu}(x)}{(q^n-q^{-n})^{2g-1}} R_{\fb}(n\nu)\in \dq\tensor
  R(T), \] where $\sum\!\!\!\!*$ means only summing finite values.
  
  Now we can apply the arguments of \secref{sec:qrat} to the sum in
  the above formula. Indeed, since $\hat R_{\fb}$ is of non-positive
  degree in $q_\nu$, the exponential convergence of the series is
  guaranteed. (There is a possible pole at 0). The residue calculations
  also go through: we apply the Residue Theorem to the form
\begin{equation}
\label{eq:def-wr}
w_{\fb}(u,x;\hbar) = \frac{e^{\{x\} u}\,du}{1-e^{u}}
\frac{\hat R_{\fb}\left(q\mapsto e^{i\pi\hbar},q_\nu\mapsto 
e^{\hbar u/2}\right)}{(e^{\hbar u/2}-e^{-\hbar u/2})^{2g-1}}.
\end{equation}
Again, the poles break into 3 parts:
\begin{enumerate}
\item $P_1=\{\ul m|\;m\in\Z,\,(q^m-q^{-m})^{1-2g}R_{\fb}(m\nu)\neq\infty\}$;
\item $P_2\subset\{\ul n/{\hbar}+\ul m|\; n\in\Z,\,n\neq0\,|m|<M\}$
  for some $M$;
\item $P_3=\{\ul m|\;m\in\Z,\,(q^m-q^{-m})^{1-2g}R_{\fb}(m\nu)=\infty\}$.
\end{enumerate}
The  set of poles $P_1$ contributes the infinite sum,  $P_2$ --
exponentially small corrections and $P_3$ -- the asymptotic
expansion. The statements in part (3) then clearly follow. \fini

\subsection{The case of several punctures.}
We will use the notational conventions from \secref{sec:3series} and
\secref{sec:modspaces}. Thus $P\subset\Sigma$ is a set of punctures
and we fix a set of conjugacy classes
$\vec\sigma:P\rightarrow\orb_\mathrm{reg}(G)$ and corresponding
elements $\vex:P\rightarrow\lt$. We set $|P|=k$, $P=\{p_1,\cdots,
p_k\}$ and use the notation $\sigma_i=\vec\sigma(p_i)$, etc. Our goal is
the study of the series \eqref{eq:trq}:
\begin{equation}
\mathrm{Tr}_q(\fb)=
\sum_{\lambda_1,..,\lambda_k \in \Omega^{+}}
H\left(\fb \prod_{i=1}^k c_{\lambda_i}^{p_i}\right)\prod_{i=1}^r 
\bar\chi_{\lambda_i}(\sigma_i).\label{eq:qtrseveralp}
\end{equation}
The analysis is similar to the one-puncture case, thus we only
highlight the differences here. The difficulty is
that this sum formally has $k$ parameters instead of 1. Our plan is to
reduce the computations to the $k=1$ case.  Again, we first state the
results which can be obtained for any $G$ and then analyze the case
$G=SU(2)$ in more detail.

We first express $H(\fb \prod_{i=1}^k \bc_{\lambda_i}^{p_i})$ in terms
of intertwiners of irreducible finite dimensional representations of
$\ue.$ Let $\Gamma$ be the exact graph with vertex $o$ and with
$2g+k-1$ edges: $\{e_i,i=1,...,2g, m_j, j=1,...,k-1\}$ 
The ribbon graph $\Gamma(\bc_\lambda^{p_i})$ is homotopic to $m_i$ in
$\pi_1(\Sigma_P)$, whereas the ribbon graph $\Gamma(\bc_\lambda^{p_k})$ is
homotopic to $e_1e_2e_1^{-1}e_2^{-1}\cdots
e_{2g-1}e_{2g}e_{2g-1}^{-1}e_{2g}^{-1}m_1\cdots m_{k-1}$. Again, we
can compute $H(\fb \prod_{i=1}^k c_{\lambda_i}^{p_i})$ using the
cutting operation.

Let $C_{\fb}(e_{2i-1})=V_{\mu(i)}, C_{\fb}(e_{2i})=V_{\nu(i)}$, let
$C_{\fb}(m_j)=V_{\zeta(j)}, j=1,...,k-1$. The element $\fb$ determines
an invariant $\phi_{\fb}\in\mathrm{Hom}(\mathcal{V},\dq)$, where
$\mathcal{V}=\tensor_{j=k-1}^1\mathcal{V}_j'\tensor \tensor_{i=g}^1
\mathcal{V}_i,$ with $\mathcal{V}_i= V_{\nu(i)}^{*}\tensor
V_{\mu(i)}^{*}\tensor V_{\nu(i)}\tensor V_{\mu(i)} $, and
$\mathcal{V}_i'= V_{\zeta(i)}^{*}\tensor V_{\zeta(i)}$.

Let $V$ be a finite dimensional $\ue$-module, let
$\eta_{V}(\lambda,\mu)$ be the elements of
$\mathrm{Hom}(V_{\lambda},V_{\lambda}\otimes V^{*}\otimes V)$
defined by $\eta_{V}(\lambda,\mu)=\sum_{j}\alpha^{j}\circ\alpha_j$
where $\sum_j \alpha^{j}\otimes \alpha_j =\delta(V;\lambda,\mu).$

We can therefore rewrite the formula for $H(\fb \prod_{i=1}^k
c_{\lambda_i}^{p_i})$ as:
\begin{multline*}
(\qdim V_\lambda)^{1-2g}\prod_{i=1}^{k-1}(\qdim V_{\lambda_i})^{-1} 
H(\fb \prod_{i=1}^k c_{\lambda_i}^{p_i})=\\ 
(\id_\lambda\otimes \phi_{\fb})
\eta_{\zeta(k-1)}(\lambda,\lambda_{k-1})\circ\cdots\circ 
\eta_{\zeta(1)}(\lambda,\lambda_{1})\circ
\xi_{\nu(g)\mu(g)}(\lambda)
\circ\cdots \circ \xi_{\nu(1)\mu(1)}(\lambda),
\end{multline*}
where we have denoted $\lambda=\lambda_k.$

The analog of \propref{thm:Hfc=Rf} holds.  Indeed if $\lambda$ is
sufficiently far from the walls of the dominant chamber, we have:
\[
\eta_{V}(\lambda, \lambda-\mu)=
\sum_{i}(\phi_{\lambda-\mu}^{v^i}\otimes \mathrm{id})
\phi_{\lambda}^{Q(\lambda)^{-1}v_i}
\]
where $\sum v^i\tensor
v_i=\delta(V[\mu])$. Denote by $W(i)$, $i=1,\cdots, 4g+2k-2$,
the appropriate weight spaces 
\begin{multline*}
  V_{\zeta(k-1)}^*[\mu_{k-1}], V_{\zeta(k-1)}[\mu_{k-1}],\cdots,
  V_{\zeta(1)}^*[\mu_{1}], V_{\zeta(1)}[\mu_{1}], V_{\nu(g)}^* [0],
  V_{\mu(g)}^*[0],
  V_{\nu(g)}[0],\\
  V_{\mu(g)}[0],\cdots, V_{\nu(1)}^*[0], V_{\mu(1)}^*[0],
  V_{\nu(1)}[0],V_{\mu(1)}[0].
\end{multline*}
 If $\lambda$ is sufficiently far from the
walls of the Weyl chamber of dominant weights, we have:
\begin{multline}
(\qdim V_\lambda)^{2g-1}\prod_{i=1}^{k-1}\qdim V_{\lambda_i} 
H\left({\bf f}\prod_{i=1}^{k-1}c^{p_i}_{\lambda-\mu_i}
c^{p_k}_{\lambda}\right)=\\
\phi_{\bf f}J_{1,\cdots,4g+2k-2}(\lambda)\prod_{j=1}^{k-1}
Q_{W(2j)}^{-1}(\lambda)\prod_{i=1}^{g}Q_{W(4i-1+2k-2)}^{-1}(\lambda)
Q_{W(4i+2k-2)}^{-1}(\lambda)\omega_{g,k}\label{eq:Hfc=Rfk-punctures}
\end{multline}
where $\omega_{g,k}$ is the diagonal element of
$\otimes_{i=1}^{4g+2k-2} W(i)$.

Again, there are appropriate rational functions $\hat R_{\fb}$ and
$R_{\fb}$, depending on the additional parameters $\mu_1,\cdots,
\mu_{k-1}$, which recover \eqref{eq:Hfc=Rfk-punctures} for
sufficiently generic $\lambda$. The proof of \propref{thm:weyl-sym}
carries over to the multi-puncture case:
\begin{prop}
\label{thm:weyl-symseveralp}
The functions $R_{\fb}(\lambda;\mu_1,\cdots,\mu_{k-1})$ satisfy 
\begin{equation}
R_{\fb}(\lambda;\mu_1,\cdots,\mu_{k-1})=
R_{\fb}(w.\lambda;w(\mu_1),\cdots,w(\mu_{k-1}))
\end{equation}
for any $w\in W_{G}$.
\end{prop}

Assume that $G=SU(2)$ in the rest of the section; $\nu$ is the
fundamental weight as before.  Denote the set of weights of
$V_{\zeta(i)}$ by $Z_i \nu$, where $Z_i$ is a finite subset of $\Z$,
invariant under $n\rightarrow -n$. The key point is that if we fix
$\lambda=n\nu$, then only finitely many possible sets of $\mu$'s can
make a non-zero contribution to \eqref{eq:qtrseveralp}. More
concretely, we need to show that
\begin{equation}
\sum_{n_j\in Z_j, n\geq N}
\frac{R_{\fb}(n\nu ;n_1 \nu,\cdots,n_{k-1}\nu)}{[n+1]_q^{2g-1}
\prod_{j=1}^{k}[n+1-n_{j}]_q}
\prod_{j=1}^{k}\frac{\sin(2\pi(n+1-n_j))x_j}{\sin(2\pi x_j)}
\end{equation}
for $N$ large enough admits an asymptotic expansion of the form stated
in \thmref{thm:main}.  This series can also be written:
\begin{equation}
\sum_{\epsilon_{j}=\pm 1, n_j\in Z_j, n\geq N}\frac{R_{\fb}((n-1)\nu ;n_1 
\nu,\cdots,n_{k-1}\nu)}{[n]_q^{2g-1}\prod_{j=1}^{k}[n-n_{j}]_q}
\prod_{j=1}^{k}\frac{\epsilon_j
e^{\epsilon_j(n-n_j)\underline{x_j}}}{2 i \sin(2\pi x_j)}. 
\end{equation}
Using the symmetry property of $R_{\fb}((n-1)\nu ;n_1
\nu,\cdots,n_{k-1}\nu),$ the series has an asymptotic expansion, if all
the series
\begin{equation}
\sum_{ n\in \Z}\csi
\frac{R_{\fb}((n-1)\nu ;n_1 \nu,\cdots,n_{k-1}\nu)}{[n]_q^{2g-1}
\prod_{j=1}^{k}[n-n_{j}]_q}
\prod_{j=1}^{k}\frac{
  e^{\sum_{j=1}^{k}\epsilon_j(n-n_j)\underline{x_j}}}
{2 i \sin(2\pi x_j)}, 
\end{equation}
have asymptotic expansions.  These series are exactly of the type
already studied in the one puncture case and they satisfy the
asymptotic property when $\sum_j \epsilon_j x_j\notin \Z$.
This last condition is exactly the condition of being ``non-special''.
\fini

\section{Conclusion}
\label{sec:conclusion}

The purpose of this section is to review what we have done so far, and
formulate the conjecture which served as the original motivation for
starting this work. In the theorem below, we put our results
together. For simplicity, we consider the one-puncture case only.
\begin{thm}
  \label{thm:final}
  Let $G=\sut$, $\Sigma$ be a compact genus $g$ Riemann surface with a
  puncture at a point $p$, and let $\sigma$ be a regular conjugacy class
  in $\sut$. 
\begin{enumerate}
\item 
Then there is a non-commutative $q$-deformation
  $A_q=F^q(\M[\Sigma_p](\sigma))$ over $\dq$ of the space of algebraic
  functions on $\M[\Sigma_p](\sigma)$. 
\item There is a cyclic functional
\[ \trq^\sigma: A_q\rightarrow \mathrm{Mer}^\mathrm{asym}_1(q) \]
with values in meromorphic functions in the variable $q$ on the unit
disc which have an asymptotic expansion in the variable $\hbar$,
related to $q$ by $q=e^{\pi i \hbar}$, as $\hbar\rightarrow i0^+$.
\item The equation
  \begin{equation}
    \label{eq:final}
c(g,q)\delta(\sigma)\ase(\trq^\sigma(1)) = {\hbar}^{g-1}\int_{\M(\sigma)}
e^{\omega_\sigma/\hbar}\hat A(\M(\sigma))
  \end{equation}
holds, where $\ase:\mathrm{Mer}^\mathrm{asym}_1\rightarrow\C\hh$ is
the asymptotic expansion, $\omega_\sigma$ is the standard
symplectic form (cf. \secref{sec:sympl-form}) and $c(g,q)$ is the
normalization constant defined in \secref{sec:qrat}.
\end{enumerate}
\end{thm}

Denote the asymptotic expansion $\ase(\trq^\sigma)$ by
$\tr_\hbar^\sigma$. To proceed we need to make the following
assumptions (cf.  \secref{sec:alg-loc}): \smallskip

1. \emph{The deformation $F^q(\M^G(\sigma))$ of the algebraic
  functions on $\M^G(\sigma)$ is local (cf. \secref{sec:alg-loc} for
  the
  definition).} \\
2.  \emph{The trace functional $s^*\tr^\sigma_\hbar$ on
  $F(\M^G(\sigma))\hh$, which is the pull-back of
  $\tr^\sigma_\hbar:\fqmgp\rightarrow\Ch$ via a local section
  $s:\fmgp\rightarrow\fqmgp\hh$ extends to all smooth functions
  $C^\infty(\M^G(\sigma))\hh$. } \smallskip

We will give a proof of these two statements in a forthcoming second
part of our paper.

This finally allows us to formulate our conjecture:
\begin{conj}
\label{thm:conj}
The characteristic class of the deformation $A_q$ is
$\omega_\sigma/\hbar$, i.e. $A_q$ is basic, and up to an appropriate
power of $\hbar$ the functional 
$f\mapsto c(g,q)\delta(\sigma)\ase(\trq^\sigma(f))$ is the canonical trace.
\end{conj}
We expect this conjecture to hold in complete generality, for all
groups and arbitrary number of punctures. Note that in order for the
asymptotic expansions to exist, one needs to assume that $\vec\sigma$
is non-special.

We do not know how to approach this conjecture at the moment. Clearly,
the work \cite{AGS} is relevant, but it is not clear how to make the
connection rigorous.  What we have shown is that the two statements in
the conjecture are consistent: \emph{If $A_q$ is basic then the
  asymptotic trace is canonical.}

Finally, note that for groups other than $\Sun$, the moduli space of
flat connections with fixed holonomies is an orbifold even for generic 
inserted conjugacy classes. Our results could shed some light on the
correct form of the orbifold version of the index theorem of Fedosov
and Nest-Tsygan.

\section{Appendix} 
\label{sec:appendix}
 
In this appendix we collected facts about quantum groups of which we
make use in the main text.  We recall the basic definitions in
\secref{sec:q-basic} and make them more explicit in the
$\mathfrak{sl}_2$ case in \secref{sec:q-sl2}. The monographs \cite{CP,K}
are the basic references for this part. In \secref{sec:q-dynqg}
we recall the definitions of dynamical quantum groups and dynamical
Weyl groups defined in \cite{EV1,EV2,EV3}. We also present some
explicit computations in the $\mathfrak{sl}_2$ case, which are based
on a brief remark in \cite{TV}, but are not available in this form in
the literature.
 
\subsection{Quantum Universal Algebra} 
 \label{sec:q-basic}
For an indeterminate $q$, define the following elements of 
$\Z[q,q^{-1}]$: $[m]_q=\frac{q^m-q^{-m}}{q-q^{-1}}$ for $m\in 
\Z$, $[n]_q !=[n]_q\cdots [1]_q$, for $n\in \N,$ and 
\[ \left[\begin{array}{c}n\\m \end{array} \right]_{q} 
=\frac{[n]_q!}{[m]_q![n-m]_q!}, 0\leq m\leq n. \] 
  
Let $\mathfrak{g}$ be a finite dimensional complex simple Lie algebra 
of rank $r$ with Cartan matrix $(a_{ij})$, and let $d_i$ be the 
coprime positive integers such that the matrix $d_i a_{ij}$ is 
symmetric. Introduce the notation $q_{i}=q^{d_i}.$ 
 
$\mathfrak{U}_q(\mathfrak{g})$ is the $\C(q)$ Hopf algebra 
generated by $K_i, K_i^{-1}, e_i, f_i,\, i=1,\dots,r$ satisfying the 
defining relations: 
\[K_i K_{i}^{-1}=K_{i}^{-1}K_i=1, K_{i}K_{j}=K_{j} K_{i},\] 
\[ K_i e_j K_{i}^{-1}=q_{i}^{a_{ij}} e_j, K_i f_j 
K_{i}^{-1}=q_{i}^{-a_{ij}} f_j, 
\]  
\[e_i f_j-f_j e_i={\delta_{ij}}\frac{K_i-K_i^{-1}}{q_{i}-q_{i}^{-1}},\] 
\[\sum_{r=0}^{1-a_{ij}}\left[\begin{array}{c}1-a_{ij}\\r \end{array} 
\right]_{q_i}(-1)^r e_{i}^{1-a_{ij}-r}e_j e_{i}^r=0,\] 
\[\sum_{r=0}^{1-a_{ij}}\left[\begin{array}{c}1-a_{ij}\\r \end{array} 
\right]_{q_i}(-1)^r f_{i}^{1-a_{ij}-r}f_j f_{i}^r=0.\] 
 
The coproduct is defined by: 
$$\Delta(e_i)=e_i\otimes 1+K_{i}^{-1}\otimes e_i,
\Delta(f_i)=f_i\otimes K_i+1\otimes f_i, \Delta(K_i)=K_i\otimes K_i.$$
We denote the counit by $\epsilon$ and the antipode by $S$. The sum of 
positive roots may be expressed as $2\rho=\sum_{i=1}^rm_i
\alpha_i$ with $m_i\in \N$. Define $K_{2\rho}=\prod_{i=1}^r
K_{i}^{m_i}$. Then for every $a\in \mathfrak{U}_q(\mathfrak{g})$
we have $S^2(a)= K_{2\rho}aK_{2\rho}^{-1}.$ 
 
As the above relations are defined over the ring ${\mathcal
  A}=\Z[q,q^{-1}]$, by adding to the list of generators the divided
powers $(K_i-K_i^{-1})/(q_i-q_i^{-1})$, one can define the
``non-restricted'' integral form of $\mathfrak{U}_q(\mathfrak{g})$, an
${\mathcal A}$-subalgebra $\mathfrak{U}_{\mathcal A}(\mathfrak{g})$ of
${\mathfrak U}_q(\mathfrak{g})$ such that the natural map
$\mathfrak{U}_{\mathcal A}({\mathfrak g})\otimes_{\mathcal
  A}\C(q)\rightarrow {\mathfrak U}_q(\mathfrak{g}) $ is an isomorphism
of $\C(q)$-algebras.  Then one can specialize $q$ to a non zero
complex number $q_0\in \C$ by $U_{q_0}(\mathfrak{g})=
\mathfrak{U}_{\mathcal A}(\mathfrak{g})\otimes_{\mathcal A} {\mathbb
  C}$ using the homomorphism $\mathrm{ev}_{q_0}:{\mathcal
  A}\rightarrow \C$ which sends $q$ to $q_0$. In particular, one can
set $q=1$ and obtain that $U_1(\lig)$ is essentially isomorphic to
$U(\lig)$. (One needs to set $K_i=1$ as well; for the details see
\cite[\S9.2]{CP}).  When we speak of a representation $V$ of the
quantum group, we assume that the operators $K_i$ act diagonalizably,
with eigenvalues $q^n$, $n\in\C$. This assures a ``good limit'' as
$q\rightarrow 1$, i.e. an appropriate action of $U(\lig)$ on $V/(q=1)$.

In the text, we enlarged the ring
$\mathcal{A}$ to the ring $\dq$, which consists of those rational
functions in $q$ which have no poles on the unit disc except possibly
at 0. To simplify the notation we denoted this algebra by
$\ue=\mathfrak{U}_{\mathcal A}(\mathfrak{g})\tensor_{\mathcal{A}}
\dq$.
 
For every complex weight $\lambda\in\mathfrak{h}$, define the Verma
module $M_\lambda$ as the universal $\ue$-module generated by a vector
$v_\lambda$ and relations 
$$K_{i}v_{\lambda}=q^{(\lambda,
  \alpha_i)}v_{\lambda},\; e_i v_{\lambda}=0 ,\quad i=1,\dots,r.$$ 
\begin{rem}
  Note that we defined $(.,.)$ to be the basic inner product on
  $\tdual$, normalized by the condition, that the long roots have
  square length 2 (cf. \secref{sec:lie}). The formulas above then work
  for simply laced Lie algebras only and we assume in what follows
  that $\lig$ is such a Lie algebra. In the non simply laced case one
  needs to normalize the inner product in such a way that the short
  roots have square length 2.
\end{rem}
If $\lambda\in \Omega^+$ is a dominant integral weight, then $M_\lambda$
has a unique finite dimensional quotient $V_{\lambda}$ which is
irreducible.

Define the $q$-dimension of $V_\lambda$ by
$\qdim(V_{\lambda})=\tr_{V_{\lambda}}(K_{2\rho})\in
\N[q,q^{-1}]$. From the classical Weyl formula we have
$\qdim(V_{\lambda})=\prod_{\alpha\in\Delta^+}
\frac{[(\lambda+\rho, \alpha)]_q}{[(\rho,\alpha)]_q}$. 

The Hopf algebra ${\mathfrak U}_q({\mathfrak g})$ has some additional
special properties: it is 
\begin{itemize}
\item \emph{quasitriangular}: there is an operator $R$ in a completion
  of ${\mathfrak U}_q({\mathfrak g})^{\otimes 2}$ such that for every
  $a\in {\mathfrak U}_q({\mathfrak g}),\; \Delta^{op}(a)=R \Delta(a)
  R^{-1}$, and $(\Delta \otimes \mathrm{id})(R)=R_{13}R_{23}$, $(
  \mathrm{id}\otimes \Delta)(R)=R_{13}R_{12}$.
\item \emph{a ribbon Hopf algebra}: it can be completed with the element
  $u=\sum_{i}S(b_i)a_i,$ and with a central element $v$ (the ribbon
  element), defined by $v^2=uS(u)$, and
  $\Delta(v)=(R_{21}R_{12})^{-1}(v\otimes v)$. The action of $v$ on
  irreducible representations $V_{\lambda}$ is the constant
  $v(\lambda)=q^{-C(\lambda)}$ where $C(\lambda)$ is the value of the
  classical quadratic Casimir in the classical representation
  associated to the dominant weight $\lambda$.  We use the symbol
  $\mu$ for the grouplike element $uv^{-1}=K_{2\rho}$. It plays an
  important role in the ``attaching the candycane'' transformation in 
  \secref{sec:cyclic}. While the natural pairing $V^*\tensor
  V\rightarrow\C(q)$ is invariant, the pairing $V\tensor V^*$ is
  not. It needs to be composed with $\mu^{-1}$ acting on the second
  factor. 
\end{itemize}
This elements act in the appropriate modules over $\ue$ as well.
 
\subsection{The example  $\mathfrak{U}_q(sl(2))$} 
\label{sec:q-sl2}
In this paragaph we write down some explicit formulas for the objects
defined above in the simplest case of $\lig=\mathfrak{sl}_2$. 

The $R$-matrix belongs to a completion of ${\mathfrak
  U}_q(\mathfrak{sl}_2)^{\otimes 2}$ and is given by
$R=R^{0}q^{h\otimes h/2}$, where $R^{0}=\exp_{q^{-2}}(q e\otimes f)$
and $K=q^{h} $. The $q$-exponential is defined by
$\exp_{\alpha}(z)=\sum_{n=0}^{+\infty}z^{n}/(\alpha;\alpha)_{n}$, where
$(a;b)_{n}=\prod_{k=0}^{n-1} (1-a b^{k})$.
 
The irreducible representations $V_m$ of $\mathfrak{U}_q(\mathfrak{sl}_2)$ are 
classified by a positive integer $m$. We can choose a basis $v^{m}_k,\; 
0\leq k\leq m$ of $V_m$, on which the action of the generators of 
$\mathfrak{U}_q(\mathfrak{sl}_2)$ is 
$$Kv^{m}_k=q^{(m-2k)}v^{m}_k, \; ev^{m}_k=[m-k+1]v^{m}_{k-1}, \;
fv^{m}_k=[k+1]v^{m}_{k+1},$$ 
where we omitted the subscript from the 
definition of the $q$-integers. 
 
A basis of the Verma module $M_{\lambda}$ is denoted $u^{\lambda}_{k}$, 
$k\in \N$, on which the action of the generators of 
$\mathfrak{U}_q(\mathfrak{sl}_2)$ is 
$$Ku^{\lambda}_k=q^{(\lambda-2k)}u^{m}_k, \;
eu^{\lambda}_k=[k][\lambda-k+1]u^{\lambda}_{k-1}, \;
fu^{\lambda}_k=u^{\lambda}_{k+1}.$$ 
  
The action of the ribbon element $v$ on $V_m$ is well defined and is
simply multiplication by the constant $q^{-m(m+2)/2}$.

One can define a quantum version of the Weyl reflection as follows.
Define an algebra automorphism $T$ of
$\mathfrak{U}_q(\mathfrak{sl}_2)$ by  
$$
T(e)=-q f, \;T(f)=-q^{-1} e,\;T(q^h)=q^{-h}.$$
The operator $T$
intertwines $\Delta$, the standard coalgebra structure of
$\mathfrak{U}_q(\mathfrak{sl}_2)$, with the opposite coalgebra
structure $\Delta^{\mathrm{op}}$, i.e $\Delta^{\mathrm{op}}T=(T\otimes
T)\Delta$. Then one can define an element $w$ such that
$T(a)=waw^{-1}$ and $\Delta(w)=R^{-1}(w\otimes w)$. Its action on
$V_m$ can easily be computed: $wv^{m}_k=(-1)^{k}q^{-m^2/4-k}
v^{m}_{m-k}$. We have $w^2=(-1)^{m}\id$ on $V_m$.

\subsection{Dynamical Quantum Groups} 
\label{sec:q-dynqg}
Here we give a short survey of the formalism of fusion matrices
introduced in \cite{EV1, EV2} and dynamical quantum Weyl groups from
\cite{TV,EV3} in a form necessary for our applications. The original
papers mainly deal with intertwiners between Verma modules instead of
irrducible finite dimensional modules. We have had lots of help from
Pavel Etingof here.

\emph{Notation}. We presume the notation of \secref{sec:lie} and
\secref{sec:proptraces}.  For any $\h$-diagonalizable representation
$V$ of the algebra $\ue$ and weight $\mu\in\hdual$ denote by $V[\mu]$
the $\mu$-weight space in $V$. Thus we have $V=\osum_\mu V[\mu]$. Note
that for a highest weight representation $U_\lambda$ with highest
weight $\lambda$, we have $U_\lambda[\mu]=0$ unless $\mu\leq\lambda$
with respect to the standard partial order (cf. \secref{sec:lie}).
Also, we will denote the weight of a vector $v$ of pure weight in a
representation $V$ by $\wt(v)$.
\begin{deff}
\label{thm:exp-val}
  Let $V$ be a finite dimensional representation of $\ue$, and
  $U_{\lambda}, U_{\mu}$ be highest weight representations generated
  by the vectors $u_\lambda$ and $u_\mu$, correspondingly. The {\em
    expectation value} $\langle\phi\rangle$ of an intertwiner $\phi$
  from $\Hom(U_{\lambda},U_{\mu}\tensor V)$ is defined by the
  equation $\phi(u_\lambda)=u_{\mu}\otimes \langle\phi\rangle+\sum_{i}
  u_i\otimes v_i$, where $\wt(u_i)<\mu$. This yields a map
  $$\langle\rangle:\Hom(U_{\lambda},U_{\mu}\tensor V)\rightarrow
  V[\lambda-\mu].$$
\end{deff}

We will be interested in two cases: when $U_{\lambda}=V_\lambda$ the
irreducible highest weight representation for a dominant integral
weight $\lambda$ and when $U_{\lambda}=M_\lambda$ the Verma module for
arbitrary $\lambda\in\hdual$. We denote the intertwiners in the second
case by
$$\tilde I(V;\lambda,\mu)=\Hom(M_\lambda, M_\mu\tensor V)\text{
  and } \tilde I^*(V;\lambda,\mu)=\Hom(M_\mu, M_\lambda\tensor
V^*).$$

\begin{deff}
\label{thm:generic}
  We will say that $\mu\in\tdual$ is \emph{not} generic with respect
  to the dominant weight $\nu\in \Omega^+$ if for some positive root
  $\alpha\in\Delta^+$ we have\footnote{We remind the reader that we
    are assuming $G$ is simply laced}
\[ 0\leq (\alpha,\lambda+\rho)\alpha \leq\nu
\]
\end{deff}
\begin{lemma}
  Let $\lambda,\mu\in \Omega^+$.
 \begin{enumerate}
 \item The natural map $\pi:\tilde
   I(V,\lambda,\mu)\rightarrow\Hom(M_\lambda, V_{\mu}\tensor V)$
   factors through a map $\tilde\pi:\tilde I(V,\lambda,\mu)\rightarrow
   I(V,\lambda,\mu)$, which is compatible with the expectation value
   maps.
 \item The expectation value map
   $\langle\rangle:I(V;\lambda)\rightarrow V[\lambda-\mu]$ is
   injective.
 \item Suppose that $V=V_\nu$ with $\nu\in \Omega^+$ and that $\mu$ is
   generic with respect to $\nu$. Then the expectation value map
   $\langle\rangle:\tilde I(V_\nu;\lambda,\mu)\rightarrow
   V_\nu[\lambda-\mu]$ has a canonical (right) inverse
   $v\mapsto\phi^v_\lambda$. In particular, the map $\langle,\rangle$
   is surjective.
  \end{enumerate}
\end{lemma}
The Lemma implies the following important statement:
\begin{prop}
\label{thm:important}
  If $\mu$ is generic with respect to $\nu$, then the map
    \[\langle\rangle:I(V_{\nu};\lambda,\mu )\rightarrow
    V_{\nu}[\lambda-\mu ]\] is an isomorphism of vector spaces.
\end{prop}
\begin{rem}
\label{thm:r-generic}
  1. One can extend the notion of genericity of $\mu$
  with respect to $V_\nu$ to an arbitrary finite dimensional
  representation $V$  by additivity, and thus conclude that
  $\langle,\rangle:I(V;\lambda,\mu)\rightarrow
  V[\lambda-\mu]$ is an isomorphism. \\
  2. We will also denote by $\phi^v_\lambda$ the map
  $\tilde\pi\phi^v_\lambda\in I(V,\lambda,\lambda-\wt(v))$ if this
  causes no confusion.\\
  3. Note that in the case of $\lig=\mathfrak{sl}_2$, there are only
  finitely many weights $\mu$ non-generic with respect to a particular
  representation $V$.
\end{rem}
\emph{Proof of the Lemma}: (1). Let $\phi\in \tilde
I(V,\lambda,\mu)$ be a non-zero intertwiner. Then the image
$\pi(\phi(M_{\lambda}))$ is a finite dimensional module of highest
weight $\lambda$, which necessarily has to be isomorphic to $V_{\lambda}$.
Clearly, then $\pi$ factors through $\tilde\pi$. \\
(2). Let $\phi\in I(V,\lambda,\mu)$ and represent the image of the
highest weight vector as $\phi(v_\lambda)=\sum_{i}x_i\otimes y_i,$
where we assume that the $x_i$s are vectors of pure weight in $V_\mu$.
Split the sum as
\[ \sum_{i}x_i\otimes y_i = \sum_{j}x_j^\maxx\otimes y_i + 
\sum_{l}x_l\otimes y_l \] where the vectors $x_j^\maxx$ are the
vectors of maximal weight among those which occur in the original sum.
It follows from the intertwiner property that any such vector
$x_j^\maxx$ has to be a singular vector, i.e. be killed by all of the
$e_i$'s. Then the statment follows since $V_\lambda$ has only one such
vector, the highest weight vector, and this implies the
statement.\\
(3). This statement is a slight generalization ($\lambda, \mu$
arbitrary) of Etingof-Styrkas \cite[Proposition 2.1.]{ES}.$\Box$

The Proposition allows us to introduce the basic objects of
\cite{EV1}: fusion matrices and the dynamical Weyl group operators.
Below we will always assume that the necessary genericity conditions
hold. In particular, we have an intertwiner $\phi^{v}_{\lambda}\in
\tilde I(V,\lambda,\lambda-\wt(v))$ such that $\langle\phi\rangle =v.$
If $V, W$ are finite dimensional modules, define an endomorphism of
$V\otimes W$, called \emph{fusion matrix} and denoted
$J_{VW}(\lambda)$, by the equation
\begin{equation}
  \label{eq:fusion-matrix}
  \langle \phi_{\lambda-\wt(w)}^{v}\circ\phi_{\lambda}^{w}\rangle=
J_{VW}(\lambda)(v\otimes w),
\end{equation}
where $v, w$ are pure weight vectors.  In fact, there exists a
universal element $J(\lambda)$ in a completion of $\mathfrak{U}^{\otimes
  2}$ such that $J(\lambda)$ is represented by $J_{VW}(\lambda)$ on
the module $V\otimes W$.  
This operator satisfies the so-called ``dynamical cocycle equation'':
\[
( \mathrm{id}\otimes \Delta)J(\lambda)J_{23}(\lambda)= 
( \Delta  \otimes \mathrm{id})J(\lambda)J_{12}(\lambda-h^{(3)}), 
\] 
where the notation $h^{(3)}$ stands for the action of $h$ on the 3rd
tensor component.

One can generalize the definition of $J_{VW}$ to tensor products with
more components. These higher fusion matrices are also induced by
universal elements which may be written as 
\begin{equation}
  \label{eq:higher-j}
  J_{1,2,\cdots,N}(\lambda)=
J_{1,2\cdots N}(\lambda)\cdots J_{N-1,N}(\lambda).
\end{equation}
where $J_{1,2\cdots N}(\lambda)=(id\otimes \Delta^{(N-1)})J(\lambda)$
and $\Delta^{(p)}:\ue\rightarrow \ue^{\otimes p}$ is the iterated
coproduct.

If $v\in W[\lambda-\mu]^*, w\in W[\lambda-\mu]$, the linear map
$\tr_{W}( \phi_{\lambda}^{v}\otimes
\mathrm{id}_{W})\phi_{\lambda}^{w}$ is an intertwiner from
$V_{\lambda}$ to itself, which is necessarily the identity times a
constant, which will be denoted by the same symbol.  This defines a
non-degenenerate pairing between $W[\lambda-\mu]^*$ and
$W[\lambda-\mu]$, and we can define an invertible endomorphism
$Q_{W}(\lambda)$ of $W[\lambda-\mu]$, such that:
\begin{equation} 
\label{eq:defq}
\tr_{W}(\phi_{\lambda}^{v}\otimes\id_{W}) 
\phi_{\lambda}^{w}=\langle v, Q_{W}(\lambda)w\rangle ,  
\end{equation} 
i.e $ \langle v,Q_{W}(\lambda)w\rangle
=\tr_{W}(J_{W^{*}W}(\lambda)(v\otimes w)).$ We can also define the
universal element $Q(\lambda)$ in a completion of $\ue$ by
$Q(\lambda)=\sum_{i}S(a_i)b_i$ where $J(\lambda)=\sum_{i}a_i\otimes
b_i$. It is easy to check that this element is represented on $W$ by
$Q_{W}(\lambda).$
 
Finally, we turn to the definition of the dynamical quantum Weyl group
was introduced in \cite{TV,EV2}. For $V$ finite dimensional,
sufficiently generic $\lambda\in \Omega^+$ and an element $w\in W_G$, there
is a canonical inclusion $M_{w.\lambda}\hookrightarrow M_\lambda$,
which induces an isomorphism between $\tilde
I(V;w.\lambda)\rightarrow\tilde I(V;\lambda)$. The expectation value
map identifies these spaces with $V[\nu]$ and $V[w(\nu)]$,
correspondingly, thus this isomorphism can be represented by an
operator $A_{V,w}(\lambda):V[\nu]\rightarrow V[w(\nu)]$. Again, this
operator is induced by a universal element $A_w(\lambda)$ in a
completion of $\ue$, such that on each finite dimensional ${\mathfrak
  U}_q(\mathfrak{g})$ module $V,$ $A_w(\lambda)$ is represented by an
endomorphism $A_{V,w}(\lambda)$. The operators $A_w(\lambda)$ satisfy
the following two relations:
\begin{eqnarray} 
&&A_{ww'}(\lambda)=A_{w}(w'.\lambda)A_{w'}(\lambda), 
\forall w,w'\in W_G, l(ww')=l(w)+l(w')\\ 
&&\Delta(A_w(\lambda))J(\lambda)=J(w.\lambda)A_{w}^{(2)} 
(\lambda)A_{w}^{(1)}(\lambda-h^{(2)}).\label{eq:AJ=JA} 
\end{eqnarray}
In particular, applying $(S\otimes id)$ followed by the algebra
multiplication one obtains the relation between $Q(\lambda)$ and
$Q(w.\lambda):$
\begin{equation} 
Q(\lambda)=S(A_{w}(\lambda-wh))Q(w.\lambda)A_w(\lambda). 
\label{eq:Q=SAQA} 
\end{equation}

\subsection{Explicit computation in the $\mathfrak{sl}_2$ case}
\label{sec:q-explcit}

Here we present some explicit computations of the objects defined in
the previous paragraph in the $\mathfrak{U}_q(\mathfrak{sl}_2)$ case.
To simplify our notation, we will identify the weight $\lambda=l\nu$
with the integer $l$ and the weight $\mathrm{wt}(v)$ with the symbol
$h$. 

Let $J(\lambda)$ the fusion matrix of
$\mathfrak{U}_q(\mathfrak{sl}_2),$ it can easily be computed using the
ABRR linear equation \cite{ABRR} which reads:
\begin{equation}
J(\lambda)(1\otimes q^{(2\lambda +1) h-h^2/2})=R^{0}_{21}(1\otimes 
q^{(2\lambda +1)h-h^2/2})J(\lambda),
\end{equation} 

The computation results in the following formula: 
\begin{equation} 
J(\lambda)=\sum_{n=0}^{+\infty}\frac{q^n(1-q^{-2})^{2n}}{(q^{-2}; 
  q^{-2})_n}(f^n\otimes e^n j_n(\lambda))\label{eq:Jofsu2} 
\end{equation}
where
$j_{n}(\lambda)=\prod_{k=0}^{n-1}\left(1-q^{2(\lambda-h-k)}\right)^{-1}$.
  
Recall that $Q(\lambda)=\sum_i S(a_i)b_i$, where
$J(\lambda)=\sum_{i}a_i\otimes b_i$.  After a straightforward
computation, we obtain from \eqref{eq:Jofsu2}, that
\[
Q(\lambda)v^{m}_k={}_2\varphi_{1}(q^{-2(m-k+1)},
q^{2k};q^{2(\lambda-m+2k)})(q^{-2})v^{m}_k,
\] where ${}_2\varphi_{1}$
is the basic hypergeometric function of base $q^{-2}$ evaluated at
$q^{-2}$. Using the Heine formula \cite{GR} we obtain
\begin{equation} 
Q(\lambda)v^{m}_k=q^{-2(m-k+1)k}\frac{(q^{2(\lambda+k+1)};q^{-2})_k} 
{(q^{2(\lambda-m+2k)};q^{-2})_k}v^{m}_k. \label{eq:Qsu2} 
\end{equation}
From this explicit expression, an easy computation implies the 
relation: 
\begin{equation} 
Q(-\lambda-2) Q(\lambda+h)=q^{\frac{h^2}{2}+h}v.\label{eq:QQ=v} 
\end{equation}

If $\lambda$ is a positive integer then
$\omega_{\lambda}=\frac{f^{\lambda+1}}{[\lambda+1]!}u^{\lambda}_0$ is
a singular vector in $M_{\lambda}$ and
$\phi_{\lambda}^{v^{m}_k}(\omega_{\lambda})$ is a singular vector in
$M_{\lambda-(m-2k)}\otimes V_{m}$.  The endomorphism
$A_{V_m}(\lambda)$ is defined by
$\phi_{\lambda}^{v^{m}_k}(\omega_{\lambda})=\omega_{\lambda-(m-2k)}\otimes
A_{V_m}(\lambda)(v^{m}_k)+\text{ terms of lower weight}.$ It can be
shown that $A_{V_m}(\lambda)$ is the value of the element $A(\lambda)$
acting on $V_m,$ where
 \begin{equation} 
A(\lambda)=v^{-1} K q^{-h^2/4}w 
Q(\lambda).\label{eq:A=wQ} 
\end{equation}
In order to show this equality one can follow the proof of \cite{EV3} or 
proceed along the lines of \cite{TV}: we can first compute 
$\phi_{\lambda}^{v^{m}_k}(u^{\lambda}_0)$, and then compute 
$\phi_{\lambda}^{v^{m}_k}(\omega^{\lambda})$ by applying 
$\frac{f^{\lambda+1}}{[\lambda+1]!}$. This is leads to 
\begin{equation} 
A_{V_m}(\lambda)(v^{m}_k)=\frac{[k+n]!}{[k]!}  
{}_3\varphi_{2}(q^{2(\lambda-n+1)}, 
q^{2(-m+k-1)},q^{2k};q^{-2n-2}, q^{2(\lambda-n)})(q^{-2})v^{m}_{m-k}, 
\end{equation} 
where $n=m-2k.$ 
 
By using the q-analog of Saalsch\"utz formula \cite{GR}, we arrive at: 
\begin{equation} 
A_{V_m}(\lambda)(v^{m}_k)=(-1)^kq^{-k(m-k+1)} 
\frac{(q^{2(\lambda+k+1)};q^{-2})_k} 
{(q^{2(\lambda-m+2k)};q^{-2})_k }v^{m}_{m-k} 
\end{equation} 
from which equation \eqref{eq:A=wQ} follows.

\bibliographystyle{unsrt}

\begin{thebibliography}{10}  
  
\bibitem{AB} M. Atiyah, R. Bott, \emph{Yang-Mills equations over Riemann
  surfaces}, Philos. Trans. Roy. Soc. London Ser. A
 \textbf{308} (1983), no. 1505, 523-615.

\bibitem{ABRR} D. Arnaudon, E. Buffenoir, E. Ragoucy, Ph. Roche,
  \emph{Universal solutions of quantum Yang-Baxter equations},
  Lett. Math. Phys. \textbf{44} (1998), no. 3, 201-214.

\bibitem{AGS} A. Y. Alekseev, H. Grosse, V. Schomerus,
  \emph{Combinatorial Quantization of the Hamiltonian Chern-Simons
    Theory I}, Comm. Math. Phys. {\bf 172} (1995) 317-358.\\
  A. Y. Alekseev, H. Grosse, V. Schomerus, \emph{Combinatorial
    Quantization of the Hamiltonian Chern-Simons Theory II,}
  \newblock{Comm. Math. Phys. {\bf 174} (1995) 561-604. }

\bibitem{AMR} J. Andersen, J. Mattes, N. Reshetikhin, \emph{The
    Poisson structure on the moduli space of flat connections and
    chord diagrams,} \newblock{Topology {\bf 35} (1996), 1069-1083.}

\bibitem{AMRq} J. Andersen, J. Mattes, N. Reshetikhin,
  \emph{Quantization of the algebra of chord diagrams},
  Math. Proc. Cambridge Philos. Soc. \textbf{124} (1998), no. 3, 451-467.

\bibitem{Be} A. Beauville, \emph{Vector bundles on curves and
    generalized theta functions: recent results and open problems},
  Math. Sci. Res. Inst. Publ.  \textbf{28}, Cambridge University
  Press, Cambridge 1995, pp. 17-33. (Preprint: alg-geom/9404001)

\bibitem{BH} H. Boden, Yi Hu, \emph{Variations of moduli of parabolic
    bundles}, Math. Ann. \textbf{301}, no. 3, 539-559.

\bibitem{BL}
J.-M. Bismut, F. Labourie,
\emph{Symplectic geometry and the Verlinde Formulas,}
\newblock{Preprint Universit\'e de Paris-Sud, (1998).}


\bibitem{BR} E. Buffenoir, Ph. Roche, \emph{Two dimensional
    lattice gauge theory based on a quantum group,} \newblock{
    Comm. Math. Phys.} {\bf 170} (1995).  \\
E. Buffenoir, Ph. Roche,
  \newblock{Link Invariants and Combinatorial Quantization of
    Hamiltonian Chern Simons Theory,} \newblock{ Comm. Math. Phys.
    {\bf 181}  (1996), 331-365.}
 
\bibitem{CP}
V. Chari, A. Pressley,
\newblock{Quantum Groups,}
\newblock{Cambridge University Press (1994).}

\bibitem{De} P. Deligne, \emph{D\'eformations de l`alg\`ebre des
    fonctions d`une vari\'et\'e symplectique: comparaison entre
    Fedosov et De Wilde, Lecomte}, Selecta Math. (N.S.) \textbf{1}
  (1995), no. 4, 667-697.

\bibitem{DL} M. De Wilde, P. Lecomte, \emph{Existence of
    star-Products and of formal deformations of the Poisson Lie
    algebra of arbitrary symplectic manifolds,} \newblock{Lett. Math.
    Phys. {\bf 7} (1983), 487-496.}

\bibitem{ES} P. Etingof, K. Styrkas, \emph{Algebraic Integrability
    of Macdonald Operators and Representation Theory of Quantum
    Groups,} \newblock{Comp. Math. {\bf 114} (1998), 125-152.}
  \newblock{Preprint: math.QA/9603022}


\bibitem{EV1}  
P. Etingof, A. Varchenko,  
\emph{Exchange Dynamical Quantum Groups,}  
\newblock{Comm. Math. Phys. {\bf 205} (1999), 19-52,}  
\newblock{Preprint: math.QA/9801135.}  

\bibitem{EV2} P. Etingof, A. Varchenko, \emph{Traces of
    Intertwining Operators for Quantum Groups and Difference Equations
    I,} \newblock{Preprint: math.QA/9907181}.

\bibitem{EV3}
P. Etingof, A. Varchenko,  
\emph{Dynamical Weyl Groups and Applications,}  
\newblock{in preparation}  
  
  
\bibitem{Fe}
B. Fedosov,
\newblock{Deformation Quantization and Index Theory,}
\newblock{Mathematical Topics, 9}
\newblock{Akademie Verlag, Berlin (1996)}
 
\bibitem{FR} V. V. Fock, A. A. Rosly, \emph{Poisson structures on
    moduli space of flat connections on Riemann surfaces and
    the $r$-matrices}, Moscow Seminar in Math. Physics, 67-86,
  Amer. Math. Soc. Transl. Ser. 2 \textbf{191}(1999),
  \newblock{q-alg/9802054}

\bibitem{G} V. Guillemin, \emph{Star products on prequantizable symplectic
  manifolds}, Lett. Math. Phys. {\bf 35} (1995).

\bibitem{GR}
G. Gasper, R. Rahman,
\newblock{Basic Hypergeometric Series,}
 \newblock{Encyclopedia of mathematics and its applications,  \textbf{35}, Cambridge University Press, (1990).}

\bibitem{JK} L. C. Jeffrey, F. C. Kirwan, Intersection theory on
  moduli spaces of holomorphic bundles of arbitrary rank on a Riemann
  surface, Ann. Math. {\bf 148} (1998).

\bibitem{KS} A. V. Karabegov, M. Schlichenmaier, \emph{Identification of
  Berezin-Toeplitz deformation quantization}, Preprint: math.QA/00006063.

\bibitem{K} C. Kassel, Quantum Groups, Springer-Verlag (1995).

\bibitem{Ko}
M. Kontsevich,
\emph{Deformation quantization of Poisson manifolds, I},
\newblock{Preprint: math.QA/9709040.}

\bibitem{KL} K. Liu,  Heat kernel and moduli spaces
  I.-II. Math. Res. Letters {\bf 3} (1996), 743-762.,
Math. Res. Letters {\bf 4} (1997), 569-588.

\bibitem{MW} E. Meinrenken, C. Woodward, \emph{Hamiltonian loop group
  actions and Verlinde factorization},
J. Diff. Geom. \textbf{50}(1998), no.3, 417-469.

\bibitem{NT}
R. Nest, B. Tsygan,
\emph{Algebraic index theorem}
\newblock{ Comm. Math.  Phys. {\bf 172} (1995), no. 2, 223--262.} 


 
\bibitem{PR}
M. Polyak, N. Reshetikhin,
\emph{On 2D Yang-Mills theory and invariants of links,}
\newblock{ Deformation theory and
Symplectic geometry (Ascona, 1996), 223--248, Math. Phys. Stud.,{\bf  20},
 Kluwer Acad. Publ., Dordrecht
}

\bibitem{RT}
N. Reshetikhin, V. Turaev,
\emph{Ribbon Graphs and their invariants derived from quantum Groups,}  
\newblock{Comm. Math. Phys. {\bf  127} (1990), 1-26.} 

\bibitem{Sz1}
A. Szenes, 
\emph{The combinatorics of the Verlinde formulas,}
\newblock{Vector Bundles in Algebraic Geometry, (Durham 1993),} 
\newblock{London Math.Soc Lecture notes Ser. {\bf 208}, Cambridge
University, Press, Cambridge (1995).}

 
\bibitem{Sz2}
A. Szenes, 
\emph{Iterated residues and multiple Bernoulli Polynomials,}
\newblock{Internat. Math. Res. Notices, no. 18, (1998). } 
 
\bibitem{Sz3} A. Szenes, \emph{A vanishing theorem for trigonometric
    sums}, in preparation.

\bibitem{Te} C. Teleman, \emph{The quantization conjecture revisited}, 
  Preprint: math.AG/9808029.
  
\bibitem{Th} M. Thaddeus, \emph{Geometric invariant theory and flips},
  J. Amer. Math. Soc. {\bf 9} (1996).

\bibitem{TV} V. Tarasov, A. Varchenko, \emph{Difference Equations
    compatible with trigonometric KZ differential equations,}
  \newblock{Preprint: math.QA/0002132}

\bibitem{Tu}
V. Turaev,
\emph{Algebras of loops on surfaces, algebras of knots, and
quantization.}
\newblock{Braid group, knot theory and statistical mechanics, II,
 324--360, Adv.  Ser. Math. Phys.,
17, World Sci. Publishing, River Edge, NJ, 1994.}

\bibitem{Ve} Verlinde, \emph{Fusion rules and moduli transformations in 2D
  conformal field theory}, Nuclear Phys. B  \textbf{300} (1988), no.3,
360-376. 

\bibitem{We}
A. Weinstein,
\newblock{Deformation Quantization,}
\newblock{S\'eminaire Bourbaki, Vol. 1993/1994. Ast\'erisque N$\circ$
227 (1995).}

\bibitem{WeHo} A. Weinstein, Ping Xu, \emph{Hochschild cohomology and
    chracteristic classes for star-products}, Geometry of differential 
  equations, 177-194, Amer. Math. Soc. Transl. Ser. 2  \textbf{186}(1998).

\bibitem{WW} E. T. Whittaker, G. N. Watson, A course of modern
  analysis, Cambridge Mathematical library, Cambridge University Press 
  (1996)

\bibitem{Wi}
E. Witten,
\emph{On Quantum Gauge Theories in Two Dimensions,}
\newblock{Comm. Math. Phys. {\bf 141} (1991), 153-209.}


\end{thebibliography}

\end{document}